\documentclass{amsart}
\usepackage{graphicx}
\begin{document}

\noindent
{\LARGE\bf On the Translative Packing Densities of Tetrahedra

\medskip
\noindent
and Cuboctahedra}

\vspace{0.5cm} \noindent{\large Chuanming Zong}

\noindent
School of Mathematical Sciences, Peking University, Beijing 100871, China

\noindent
E-mail: cmzong@math.pku.edu.cn.

\vspace{0.5cm}
\noindent
{\bf Abstract.} In 1900, as a part of his 18th problem, Hilbert asked the question to determine the density
of the densest tetrahedron packings. However, up to now no mathematician knows the density $\delta^t(T)$ of the densest translative tetrahedron packings and the density $\delta^c(T)$ of the densest congruent tetrahedron packings. This paper presents a local method to estimate the density of the densest translative packings of a general convex solid. As an application, we obtain the upper bound in
$$0.3673469\cdots \le \delta^t(T)\le 0.3840610\cdots,$$
where the lower bound was established by Groemer in 1962, which corrected a mistake of Minkowski.
For the density $\delta^t(C)$ of the densest translative cuboctahedron packings, we obtain the upper bound in
$$0.9183673\cdots \le \delta^t(C)\le 0.9601527\cdots.$$
In both cases we conjecture the lower bounds  to be the correct answer.

\vspace{0.9cm}\noindent
{\large\bf 1. Introduction}

\vspace{0.4cm}\noindent
In the extended version of his talk presented at the ICM 1900 in Paris, Hilbert \cite{hilb01}
proposed 23 unsolved mathematical problems. At the end of his 18th problem, he asked
\lq\lq {\it How can one arrange most densely in space an infinite number of equal solids of given
form, e.g., spheres with given radii or regular tetrahedra with given edges $($or in prescribed
position$)$, that is, how can one so fit them together that the ratio of the filled to the
unfilled space may be as great as possible}?"

Hilbert may have asked these questions because packings of spheres and regular tetrahedra both have a long history.
Concerning the densest packing of spheres, in 1611 Kepler asserted that the densest sphere packing was given
by the face-centered cubic lattice packing. This assertion, known as Kepler's Conjecture, was proved by Hales
with Ferguson (see Hales \cite{hale05}, \cite{hale06} and Lagarias \cite{laga11}). Concerning the densest
packing of regular tetrahedra, Aristotle stated that regular tetrahedra fill space. This is not the case.
However, it took eighteen hundred years for Aristotle's error to be resolved. For the long history of investigation of this question, see Struik \cite{stru25} and Lagarias and Zong \cite{zong12}.

\medskip
Let $K$ denote a convex body in the three-dimensional Euclidean space $\mathbb{E}^3$, with boundary $\partial (K)$, interior ${\rm int}(K)$ and volume ${\rm vol}(K)$. In particular, let $T$, $O$, $C$ and $S$ denote a regular tetrahedron, a regular octahedron, a regular cuboctahedron and a unit sphere, respectively. Let $\delta^c(K)$, $\delta^t(K)$ and $\delta^l(K)$ denote the densities of the densest congruent packings, the densest translative packings and the densest lattice packings of $K$, respectively. For the detailed definitions, basic results and open problems about these densities we refer to \cite{bras05}, \cite{feje93}, \cite{feje71}, \cite{grub07}, \cite{grub87} and \cite{roge64}. It follows from their definitions that
$$\delta^l(K)\le \delta^t(K)\le \delta^c(K)\le 1\eqno (1.1)$$
holds for every convex body $K$. Moreover, both $\delta^l(K)$ and $\delta^t(K)$ are invariants under nonsingular affine linear transformations, while $\delta^c(K)$ for some $K$ is not. Then, Hilbert's problem can be restated as:
{\it To determine the values of $\delta^c(K)$, $\delta^t(K)$ and $\delta^l(K)$ for a given convex body $K$, such as a sphere or a regular tetrahedron.}

The first approach to Hilbert's problem was made by Minkowski \cite{mink04} in 1904. He proposed a two step program
to determine the values of $\delta^l(K)$. First, he defined
$$D(K)=\{ {\bf x}-{\bf y}:\ {\bf x}, {\bf y}\in K\}$$
and discovered that
$$(K+{\bf x})\cap (K+{\bf y})\not= \emptyset $$
{\it if and only if}
$$\left( \mbox{$1\over 2$}D(K)+{\bf x}\right)\cap \left( \mbox{$1\over 2$}D(K)+{\bf y}\right)
\not= \emptyset .$$
Therefore, for a discrete set $X$ in $\mathbb{E}^3$, $K+X$ is a packing if and only if ${1\over 2}D(K)+X$
is a packing. Consequently, he proved
$$\delta^t(K)={{2^3{\rm vol}(K)}\over {{\rm vol}(D(K))}}\cdot \delta^t(D(K))\eqno (1.2)$$
and
$$\delta^l(K)={{2^3{\rm vol}(K)}\over {{\rm vol}(D(K))}}\cdot \delta^l(D(K)).\eqno (1.3)$$
Usually, $D(K)$ is called the difference set of $K$. Clearly $D(K)$ is centrally symmetric, convex and centered at the origin. Second, when $K$ is centrally symmetric and centered at the origin, he discovered the following criterion for its densest lattice packings: {\it If $K+\Lambda $ is a lattice packing of maximal density, then $\Lambda $ has a basis $\{{\bf a}_1, {\bf a}_2, {\bf a}_3\}$ such that either
$$\{ {\bf a}_1, {\bf a}_2, {\bf a}_3, {\bf a}_1-{\bf a}_2, {\bf a}_2-{\bf a}_3, {\bf a}_3-{\bf a}_1\}\subset
\partial (2K)$$
or}
$$\{ {\bf a}_1, {\bf a}_2, {\bf a}_3, {\bf a}_1+{\bf a}_2, {\bf a}_2+{\bf a}_3, {\bf a}_3+{\bf a}_1\}\subset
\partial (2K).$$
As an application, he determined the density of the densest lattice packings of an octahedron $O$. In
other words, he proved
$$\delta^l(O)={{18}\over {19}}.\eqno (1.4)$$

On page 312 of \cite{mink04}, Minkowski wrote \lq\lq If $K$ is a tetrahedron, then ${1\over 2}D(K)$ is an
octahedron with faces parallel to the faces of the tetrahedron."  By routine computations, one can get
${\rm vol}(T)=\sqrt{2}/{12}$ and ${\rm vol}(O)=\sqrt{2}/3$ when both $T$ and $O$ have unit edges.
Then, by (1.3) and (1.4) Minkowski \cite{mink04} made a conclusion that
$$\delta^l(T)={9\over {38}}.\eqno (1.5)$$

\vspace{-0.1cm}
\begin{figure}[hc]
\includegraphics[height=4.5cm,width=10cm,angle=0]{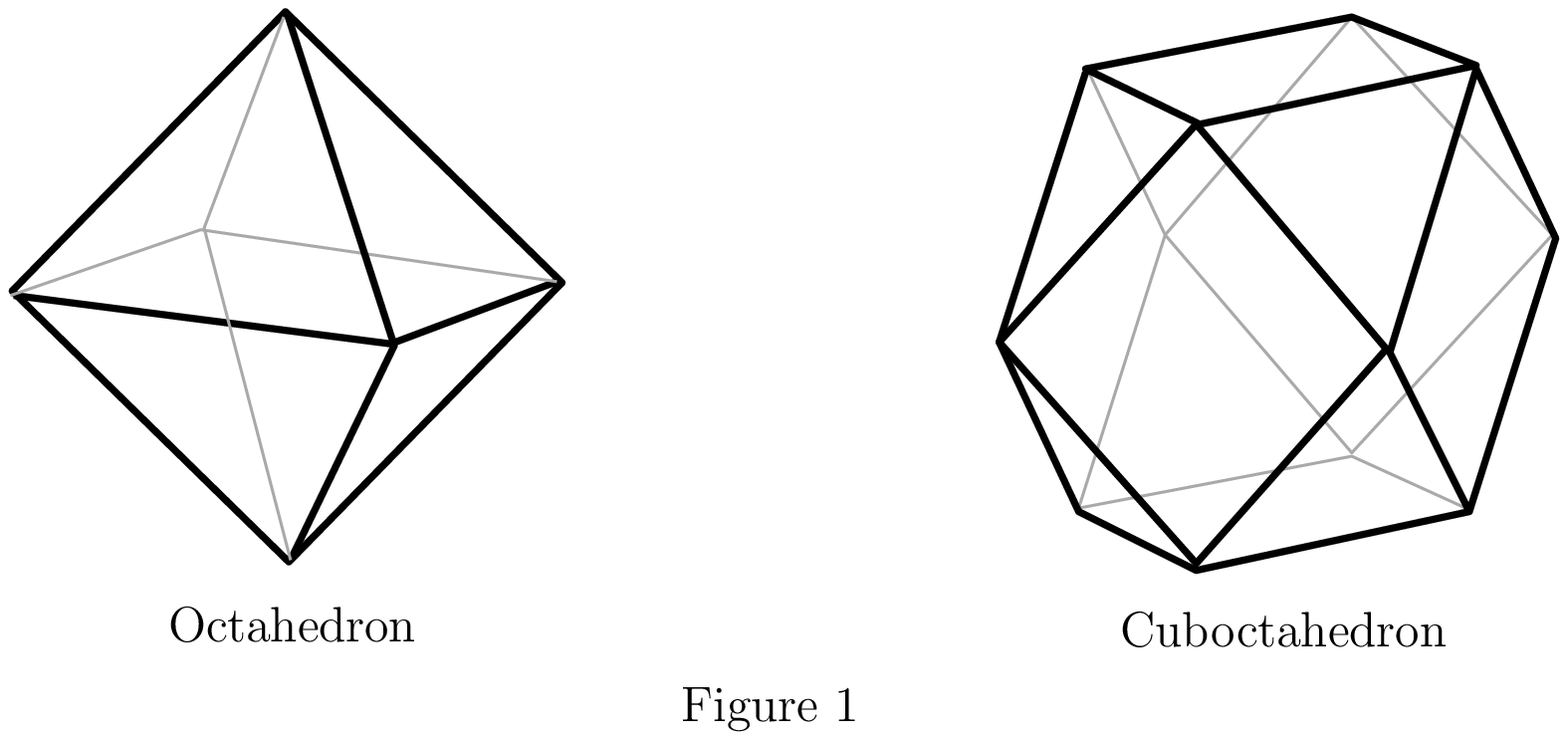}
\end{figure}

Unfortunately, in this calculation Minkowski made a mistake, which was uncovered by Groemer \cite{groe62}
in 1962. The difference set of a tetrahedron is not an octahedron, but a cuboctahedron. As is shown in Figure 1, a cuboctahedron is very different from an octahedron. In fact, it was already known to Estermann \cite{este28} and S\"uss \cite{suss28} in 1928 that
$${{{\rm vol}(D(T))}\over {{\rm vol}(T)}}=20.\eqno (1.6)$$

In 1970, Hoylman \cite{hoyl70} applied Minkowski's criterion to a cuboctahedron $C$. By considering
$38$ cases with respect to the possible positions of the three vectors of the bases, he proved
$$\delta^l (C)={{45}\over {49}}, \eqno (1.7)$$
$$\delta^l(T)={{18}\over {49}},\eqno (1.8)$$
and the optimal lattice is unique up to certain equivalence. It is noteworthy that in the densest lattice tetrahedron packing of density $18/49=0.3673\cdots $ each tetrahedron touches $14$ others; However, according to Zong \cite{zong96}, the density of the lattice tetrahedron packing of maximal kissing number $18$ is only $1/3=0.3333\cdots$.

Based on Minkowski's work, in 2000 Betke and Henk \cite{betk00} developed an algorithm by which one can determine the density of the densest lattice packing of an arbitrary three-dimensional polytope. They applied their program to recheck Hoylman's result.

In 2006, Conway and Torquato \cite{conw06} made a breakthrough in constructing dense congruent tetrahedron packings. Their idea is simple but very efficient. First, pack twenty regular tetrahedra into an icosahedron. The fraction of the icosahedral volume occupied by the tetrahedra can be $0.8567627\cdots .$ Then, construct a lattice icosahedron packing of maximum density. According to Betke and Henk \cite{betk00} it is $0.8363574\cdots .$ Thus we obtain a congruent tetrahedron packing of density approximately
$0.8363574\times 0.8567627 \approx 0.716559 .$ In other words, we have
$$\delta^c(T)\ge 0.716559 \cdots .\eqno (1.9)$$

It was conjectured by S. Ulam (see page 135 of \cite{gard01}) that {\it the maximal density $\pi /\sqrt{18}=
0.74048\cdots $ for packing congruent spheres is smaller than that for any other convex body.} Of course, it makes sense to consider regular tetrahedron as a candidate of counterexample, as Conway and Torquato \cite{conw06} did. In 2008, by constructing a cluster of eighteen congruent tetrahedra and a suitable lattice packing of the cluster,
Conway and Torquato's lower bound (1.9) was improved by Chen \cite{chen08} to
$$\delta^c(T)\ge 0.778615\cdots ,$$
which turns to support Ulam's conjecture.

Packings of regular tetrahedra may provide useful models in material science, information theory and etc. Therefore, recently it becomes an active research topic involving both mathematicians and scientists in other fields. Chen's lower bound was further improved by \cite{torq09}, \cite{torq09'}, \cite{haji09}, \cite{kall10}, \cite{chen10} and etc. So far the best known lower bound is
$$\delta^c(T) \ge {{4000}\over {4671}}= 0.856347\cdots , $$
which was discovered in 2010 by Chen, Engel and Glotzer \cite{chen10}.

On the other hand, the situation about upper bounds for $\delta^c(T)$ is rather embarrassing. It was uncovered in the fifteenth century that regular tetrahedra can not fill the space, which means that Aristotle's assertion is wrong (see \cite{zong12}). In 1961, Schmidt \cite{schm61} proved that $K$ fills the space if and only if $\delta^c(K)=1$, which implies
$$\delta^c(T)<1.$$
However, up to now the best known upper bound is only
$$\delta^c(T)\le 1-2.6\times 10^{-25},$$
which was achieved by Kallus, Gravel and Elser \cite{kall10} in 2010.

Perhaps, to determine the value of $\delta^t(K)$ is not as challenging as that for $\delta^c(K)$. However, since
$\delta^t(K)$ is invariant under nonsingular affine linear transformations, it is important. For $\delta^t(T)$ and
$\delta^t(C)$, by (1.1), (1.2), (1.6), (1.7) and (1.8) one can deduce
$${{45}\over {49}}\le \delta^t(C)\le 1\eqno (1.10)$$
and
$${{18}\over {49}}\le \delta^t(T)\le {2\over 5}.\eqno (1.11)$$
These facts support the conjecture that

$$\delta^t(C)={{45}\over {49}} \qquad {\rm and}\qquad \delta^t(T)={{18}\over {49}}.$$

\bigskip
In this paper, we present a local method to estimate the value of $\delta^t(K)$ for a general convex body $K$. First, for each translate $K+{\bf x}_i$ of a packing we define a shadow region with respect to a given direction ${\bf v}$. Second, by minimizing the volume of this shadow or an average over several particular directions, we obtain measures for the gaps of the packing. Then, upper bounds for the packing density $\delta^t(K)$ can be deduced from ${\rm vol} (K)$ and these measures. As an application to $\delta^t(C)$ and consequently to $\delta^t(T)$, we prove the following result:

\medskip\noindent
{\bf The Main Theorem.}
$$\delta^t(C)\le {{90\sqrt{10}}\over {95\sqrt{10}-4}} \qquad {\it and }\qquad \delta^t(T)\le {{36\sqrt{10}}\over {95\sqrt{10}-4}}.$$

\bigskip
Combined with (1.10) and (1.11) we get

\medskip
\noindent
{\bf Corollary 1.1.}
$$0.9183673\cdots \le \delta^t(C)\le 0.9601527\cdots $$
{\it and}
$$0.3673469\cdots \le \delta^t(T)\le 0.3840610\cdots .$$

\medskip\noindent
{\bf Remark 1.1.} {\it To read this paper, several cuboctahedron models can be helpful.}

\vspace{0.9cm}\noindent
{\large\bf 2. Methodology and Terminology}

\vspace{0.4cm}\noindent
In this section we introduce a local method to achieve upper bounds for $\delta^t(K)$. For convenience, let $\{ {\bf e}_1, {\bf e}_2, {\bf e}_3\}$ be an orthonormal basis of $\mathbb{E}^3$,
and write
$$S=\{ (x, y, z):\ x^2+y^2+z^2\le 1\}$$
and
$$W=\left\{ (x, y, z):\ |x| \le 1,\ |y| \le 1,\ |z| \le 1\right\}.$$
According to John's theorem (see page 13 of \cite{grub87}), for every three-dimensional convex body $K$ there is a nonsingular affine linear transformation $\sigma $ from $\mathbb{E}^3$ to $\mathbb{E}^3$ such that
$$\mbox{$1\over 3$}S\subseteq \sigma (K)\subseteq S$$
and therefore
$$\mbox{${\sqrt{3}}\over 9$}W\subseteq \sigma (K)\subseteq W.$$
On the other hand, it is well-known that
$$\delta^t(\sigma (K))=\delta^t(K)$$
holds for every convex body $K$ and any nonsingular affine linear transformation $\sigma $. Thus, to study $\delta^t(K)$, it is sufficient to work on convex bodies $K$ satisfying
$$\mbox{${\sqrt{3}}\over 9$}W\subseteq K\subseteq W.\eqno (2.1)$$

In particular, we define
$$C=\left\{ (x, y, z):\ \max \{|x|, |y|, |z|\}\le 1,\ | x|+|y|+|z|\le 2\right\}.$$
In fact, the cuboctahedron $C$ can be obtained from the cube $W$ by cutting off eight orthogonal unit tetrahedra. Thus we have
$${\rm vol}(C)=8\left(1-{1\over 6}\right)={{20}\over 3}.\eqno (2.2)$$

Take $K$ to be a convex body satisfying (2.1) and assume that $X=\{ {\bf x}_0, {\bf x}_1, {\bf x}_2, \cdots \}$ is a discrete set of points such that $K+X$ is a translative packing in $\mathbb{E}^3$ and let $\mathcal{X}$ denote the family of all such sets. For convenience we take ${\bf x}_0$ to be the origin ${\bf o}$ of the space and denote
the upper density of $K+X$ by $\delta (K,X)$. In other words,
$$\delta (K, X)=\limsup_{\ell\to\infty} {{n(\ell )\cdot {\rm vol}(K)}\over {{\rm vol}(\ell W)}},\eqno (2.3)$$
where $\ell $ is a positive number and $n(\ell )$ is the number of the points in $X\cap \ell W.$ It is known that
$$\delta^t(K)=\sup_{\mathcal{X}}\left\{\delta (K, X):\ X\in \mathcal{X}\right\}.\eqno (2.4)$$

Let ${\bf v}$ be a unit vector and let $s(K,X,{\bf v},{\bf x})$ denote the set of points ${\bf y}$ such that ${\bf y}={\bf x}+\tau {\bf v}$ holds for some positive number $\tau $ and the whole open segment $({\bf x}, {\bf y})$ belongs to $\mathbb{E}^3\setminus \{ K+X\}$. Then, we define the shadow region $D(K,X,{\bf v},{\bf x}_j)$ of $K+{\bf x}_j$ in the direction ${\bf v}$ to be the closure of
$$\biggl(\bigcup_{{\bf x}\in \partial (K+{\bf x}_j)}s(K,X,{\bf v}, {\bf x})\biggr) \cap (2W+{\bf x}_j).$$
The cube $2W+{\bf x}_j$ here is a localizer, which simplifies the computation. Clearly $D(K,X,{\bf v}, {\bf x}_j)$ is a measurable set associated to $K+{\bf x}_j$ and, for fixed ${\bf v}$,
$${\rm int}(D(K,X,{\bf v},{\bf x}_j))\cap {\rm int}(D(K,X,{\bf v},{\bf x}_k))=\emptyset $$
holds for any pair of distinct indices $j$ and $k$. Therefore the fraction
$${{{\rm vol}(K)}\over {{\rm vol}(K)+{\rm vol}(D(K,X,{\bf v},{\bf x}_j))}}$$
measures the local density of $K+X$ at $K+{\bf x}_j$.

\medskip
\noindent
{\bf Remark 2.1.} {\it Similar idea for two-dimensional packings, without the localizer, can be traced back to L. Fejes T\'oth \cite{feje83}.}

\medskip\noindent
{\bf Example 2.1.}  {\it When ${\bf a}_1=(2,-{1\over 3}, -{1\over 3}),$ ${\bf a}_2=(-{1\over 3}, 2, -{1\over 3})$,
${\bf a}_3=(-{1\over 3}, -{1\over 3}, 2)$ and
$$\Lambda_1=\left\{\sum z_i{\bf a}_i:\ z_i\in \mathbb{Z}\right\},$$
we get
$${\rm vol}(D(C,\Lambda_1, {\bf e}_1, {\bf u}_j))={{16}\over {27}}$$
and}
$$\delta (C,\Lambda_1)={{{\rm vol}(C)}\over {{\rm vol}(C)+{\rm vol}(D(C,\Lambda_1,{\bf e}_1,{\bf o}))}}={{45}\over {49}}.$$

\bigskip
Now we introduce and study the shadow measure $\mu (K, {\bf v})$. We define
$$\mu (K, X, {\bf v})=\min_{{\bf x}_j\in X} {\rm vol} (D(K,X, {\bf v}, {\bf x}_j))\eqno (2.5)$$
and
$$\mu (K, {\bf v})=\min_{X\in \mathcal{X}} \mu (K, X, {\bf v}).\eqno (2.6)$$
By (2.3), (2.4) and the definition of $D(K, X, {\bf v}, {\bf x}_j)$, we get
$$\delta (K, X)=\limsup_{\ell\to\infty} {{n(\ell )\cdot {\rm vol}(K)}\over {{\rm vol}((\ell +1) W)}}
\le {{{\rm vol}(K)}\over {{\rm vol} (K)+\mu (K,X,{\bf v})}}$$
and therefore
$$\delta^t(K)=\sup_{\mathcal{X}}\{ \delta (K, X):\ X\in \mathcal{X}\}\le
{{{\rm vol}(K)}\over {{\rm vol} (K)+\mu (K, {\bf v})}}.\eqno (2.7)$$

\medskip
Thus, we have proved the following result, which is one of the keys for this paper.

\medskip
\noindent
{\bf Lemma 2.1.} {\it For every convex body $K$ and every unit vector ${\bf v}$ we have}
$$\delta^t(K)\le {{{\rm vol}(K)}\over {{\rm vol}(K)+\mu (K, {\bf v})}}.$$

\medskip
It is easy to see that, for a fixed convex body $K$, $\mu (K, {\bf v})$ is a continuous function of ${\bf v}$ restricted on $\partial (S)$. If one can find a direction ${\bf v}$ with larger $\mu (K,{\bf v})$, he will be able to achieve a better upper bound for $\delta^t(K)$. Nevertheless, to determine or estimate the value of $\mu (K, {\bf v})$ itself is a hard job.

\medskip
\noindent
{\bf Lemma 2.2.} {\it For any pair of a fixed three-dimensional convex body $K$ satisfying $(2.1)$ and a fixed unit vector ${\bf v}$ there is a suitable finite discrete set $X_1$ satisfying}
$$\mu (K, {\bf v})={\rm vol}(D(K,X_1,{\bf v},{\bf o})).$$

\medskip\noindent
{\bf Proof.} Suppose that $Q_1$, $Q_2$, $Q_3$, $\cdots$ is a sequence of discrete sets in $\mathbb{E}^3$ and
${\bf q}_1$, ${\bf q}_2$, ${\bf q}_3$, $\cdots$ is a corresponding sequence of points such that $K+Q_i$ are
packings in $\mathbb{E}^3$, ${\bf q}_i\in Q_i$, and
$$\lim_{i\to\infty } {\rm vol}(D(K,Q_i, {\bf v}, {\bf q}_i))=\mu (K, {\bf v}).\eqno (2.8)$$
By defining
$$Q'_i=\left \{ {\bf q}-{\bf q}_i:\ {\bf q}\in Q_i\right\}$$
and
$$Q^*_i=Q'_i\cap 3W,$$
it follows from the definition of $D(K,X,{\bf v}, {\bf x}_j)$ that
$$D(K,Q^*_i, {\bf v}, {\bf o})=D(K,Q'_i, {\bf v}, {\bf o})=D(K,Q_i,{\bf v},{\bf q}_i)-{\bf q}_i$$
and therefore by (2.8)
\begin{eqnarray*}
\lim_{i\to\infty} {\rm vol}(D(K,Q^*_i, {\bf v}, {\bf o}))&\hspace{-0.2cm}=&\hspace{-0.2cm}\lim_{i\to\infty } {\rm vol}(D(K,Q'_i, {\bf v}, {\bf o}))\\ &\hspace{-0.2cm}=&\hspace{-0.2cm}\lim_{i\to\infty } {\rm vol}(D(K,Q_i, {\bf v}, {\bf q}_i))\\
&\hspace{-0.2cm}=&\hspace{-0.2cm}\mu (K, {\bf v}).
\end{eqnarray*}

Notice that
$$K+{\bf q}^*\subset 4W$$
whenever ${\bf q}^*\in Q^*_i$. Let $|X|$ denote the number of the points of $X$, by (2.1) we obtain
$$|Q^*_i|\le {{{\rm vol}(4W)}\over {{\rm vol}(K)}}= {{8^3\cdot 9^3}\over {\sqrt{3}^3}}\le 71,832.$$
Therefore, by a suitable selection process we can obtain a subsequence $Q^*_{i_1}$, $Q^*_{i_2}$, $Q^*_{i_3}$, $\cdots$ of the sequence $Q^*_1$, $Q^*_2$, $Q^*_3$, $\cdots$ and a certain discrete set $X_1$ such that
$$\lim_{j\to \infty }Q^*_{i_j}=X_1.$$
Clearly, $K+X_1$ is a packing in $\mathbb{E}^3$ and
$${\rm vol}(D(K,X_1, {\bf v}, {\bf o}))=\lim_{j\to\infty} {\rm vol}(D(K,Q^*_{i_j}, {\bf v}, {\bf o}))=\mu (K, {\bf v}).$$

The lemma is proved. \hfill{$\Box$}

\medskip
Lemma 2.1 and Lemma 2.2 provide a mean to determine or estimate the values of $\mu (K, {\bf v})$ and $\delta^t(K)$. Perhaps, for some particular convex body $K$ and corresponding vector ${\bf v}$, it happens that
$$\delta^t(K)={{{\rm vol}(K)}\over {{\rm vol}(K)+\mu (K,{\bf v})}}.$$
However, for most $K$, we are not so lucky. Based on (2.5), for a set $V=\{ {\bf v}_1, {\bf v}_2, \cdots, {\bf v}_n\}$ of $n$ unit vectors, we define
$$\overline{\mu }(K, X, V)=\min_{{\bf x}_j\in X}{1\over n}\sum_{i=1}^n {\rm vol} (D(K,X, {\bf v}_i, {\bf x}_j))\eqno (2.9)$$
and
$$\overline{\mu }(K, V)=\min_{X\in \mathcal{X}} \overline{\mu }(K, X, V).\eqno (2.10)$$
Similar to Lemma 2.1 and Lemma 2.2, one can prove the following generalizations.

\medskip
\noindent
{\bf Lemma 2.1$^*$.} {\it For every convex body $K$ and every set $V$ of unit vectors we have}
$$\delta^t(K)\le {{{\rm vol}(K)}\over {{\rm vol}(K)+\overline{\mu }(K, V)}}.$$

\medskip
\noindent
{\bf Lemma 2.2$^*$.} {\it For any pair of a fixed three-dimensional convex body $K$ satisfying $(2.1)$ and a fixed set $V=\{ {\bf v}_1, {\bf v}_2, \cdots , {\bf v}_n\}$ of $n$ unit vectors there is a suitable finite discrete set $X_2$ satisfying}
$$\overline{\mu }(K, V)={1\over n}\sum_{i=1}^n{\rm vol}(D(K,X_2,{\bf v}_i,{\bf o})).$$

\medskip
\noindent
{\bf Example 2.2.} {\it If we take ${\bf b}_1=(2,0,0)$, ${\bf b}_2=(1, {3\over 2},
{3\over 2})$, ${\bf b}_3=(1,-{3\over 2}, {3\over 2})$, $V=\left\{ {\bf e}_1, {\bf e}_2, {\bf e}_3\right\}$, and
$$\Lambda_2=\left\{ \sum z_i{\bf b}_i:\ z_i\in \mathbb{Z}\right\},\eqno (2.11)$$
then we have
$$\mu (C,\Lambda_2, {\bf e}_1)={1\over 3}$$
and}
$$\overline{\mu}(C,\Lambda_2, V)={1\over 3}\left({1\over 3}+2+{1\over 4}+2+{1\over 4}\right)={{29}\over {18}}.$$

\medskip
Perhaps, to study $\overline{\mu}(K,V)$ is not only a way to estimate $\delta^t(K)$, but also can determine it.
Of course, to determine the value of $\overline{\mu}(K,V)$ itself is very challenging. Nevertheless, we have the following conjecture.

\medskip
\noindent
{\bf Conjecture 2.1.} {\it There is a positive integer $\kappa$ $($possiblely $\kappa =3$$)$ such that, for each
three-dimensional centrally symmetric convex body $K$ there is a corresponding set $V$ of unit vectors satisfying both $|V|\le \kappa $ and}
$$\delta^t(K)={{{\rm vol}(K)}\over {{\rm vol}(K)+\overline{\mu}(K,V)}}.$$

\bigskip
In the rest of this paper, we will deal with the particular case $K=C$. Now we present two simple observations which will be frequently used in checking certain polytope is in the gaps of $\{ C+X\}$. Assume that $X=\{ {\bf x}_0, {\bf x}_1, {\bf x}_2, \cdots \}$ is a discrete set such that $C+X$ is a cuboctahedron packing in $\mathbb{E}^3$ and let $F$ be a certain set of points. For convenience, we use ${\bf x}_i{\bf x}_j\prec F$ to abbreviate the statement that ${\bf x}_i+\lambda ({\bf x}_j-{\bf x}_i)\in F$ holds for some positive number $\lambda $.

\medskip\noindent
{\bf Lemma 2.3.} {\it If ${\bf x}_i{\bf x}_j\prec F_i$, where $F_i$ is a facet of $C+{\bf x}_i$, then the hyperplane generated by $F_i$ separates ${\rm int}(C)+{\bf x}_i$ and ${\rm int}(C)+{\bf x}_j$.}

\medskip
The fact is obvious. A proof is not necessary.

\medskip\noindent
{\bf Lemma 2.4.} {\it If $F$ is a triangular facet of $C$, two interiorly disjoint translates of $C$ can not simultaneously touch $C$ at the interior of $F$.}

\medskip\noindent
{\bf Proof.} Without loss of generality, we assume that $F$ is the facet with vertices $(1,0,1)$, $(1,1,0)$ and $(0,1,1)$. Let $D$ denote the hexagon with vertices $(-1,0,1)$, $(0,-1,1)$, $(1,-1,0)$, $(1,0,-1)$, $(0,1,-1)$
and $(-1,1,0)$, and let $\rho ({\bf x}, {\bf y})$ denote the metric defined by $D$ on any hyperplane which is parallel with $D$. Clearly $F$ and $D$ are parallel to each other. If two translates $C+{\bf x}_i$ and $C+{\bf x}_j$ can simultaneously touch $C$ at the interior of $F$. Then, both ${1\over 2}{\bf x}_i$ and ${1\over 2}{\bf x}_j$ are interior points of $F$, and ${\bf x}_i{\bf x}_j$ is parallel with the hyperplane generated by $D$. Thus we get
$$\rho (\mbox{$1\over 2$}{\bf x}_i, \mbox{$1\over 2$}{\bf x}_j)<1,$$
$$\rho ({\bf x}_i, {\bf x}_j)<2$$
and therefore
$$({\rm int}(C)+{\bf x}_i)\cap ({\rm int}(C)+{\bf x}_j)\not= \emptyset ,$$
which contradicts the assumption of the lemma. Lemma 2.4 is proved.\hfill{$\Box$}

\medskip
Let us end this section by a conjecture and a corresponding remark.

\medskip\noindent
{\bf Conjecture 2.2.} {\it For $V=\{ {\bf e}_1, {\bf e}_2, {\bf e}_3\}$, we have}
$$\overline{\mu }(C, V)={{16}\over {27}}.$$

\medskip
\noindent
{\bf Remark 2.2.} {\it This conjecture implies both}
$$\delta^t(C)={{45}\over {49}}\qquad and \qquad \delta^t(T)={{18}\over {49}}.$$

\vspace{0.9cm}\noindent
{\large\bf 3. Cuboctahedral Packings, An Observation}

\vspace{0.4cm}\noindent
Let us divide the set $R=\{ (x,y,z)\in \partial (C):\ x\ge 0\}$ into four parts
$$R_1=\left\{ (x,y,z)\in R:\ y\ge 0,\ z\ge 0\right\},$$
$$R_2=\left\{ (x,y,z)\in R:\ y\ge 0,\ z\le 0\right\},$$
$$R_3=\left\{ (x,y,z)\in R:\ y\le 0,\ z\le 0\right\},$$
and
$$R_4=\left\{ (x,y,z)\in R:\ y\le 0,\ z\ge 0\right\}.$$
Recall that $s(C,X,{\bf e}_i,{\bf x})$ is the longest open segment $({\bf x}, {\bf y})$ such that ${\bf y}={\bf x}+\tau {\bf e}_i$ holds with some $\tau >0$ and
$$({\bf x}, {\bf y})\subset \mathbb{E}^3\setminus \{ C+X\}.$$
Then, for $i=1, 2, 3$ and $4$,  we define $D_i(C,X)$ to be the closure of
$$\biggl(\bigcup_{{\bf x}\in R_i}s(C,X,{\bf e}_1, {\bf x})\biggr)\cap 2W$$
and define
$$\mu_i(C)=\min_{X\in \mathcal{X}} {\rm vol}(D_i(C,X)).\eqno (3.1)$$

\medskip
Similar to Lemma 2.2, one can deduce the following fact:

\medskip
\noindent
{\bf Lemma 3.1.} {\it For $i=1, 2, 3$ and $4,$ there are suitable finite discrete sets $X_i$ satisfying}
$$\mu_i(C)={\rm vol} (D_i(C,X_i)).$$

\medskip
It is easy to see that
$$D(C,X,{\bf e}_1, {\bf o})=\bigcup_{i=1}^4D_i(C,X)\eqno (3.2)$$
and
$${\rm int}(D_i(C,X))\cap {\rm int}(D_j(C,X))=\emptyset \eqno (3.3)$$
holds for distinct $i$ and $j$. On the other hand, by symmetry, we have
$$\mu_1(C)=\mu_2(C)=\mu_3(C)=\mu_4(C).\eqno (3.4)$$
Therefore, by (2.6), (2.5), (3.1)-(3.4) we get
\begin{eqnarray*}
\hspace{5.3cm}\mu (C, {\bf e}_1)&\hspace{-0.2cm}=&\hspace{-0.2cm}\min_{X\in \mathcal{X}} {\rm vol}(D(C,X,{\bf e}_1, {\bf o}))\\
&\hspace{-0.2cm}=&\hspace{-0.2cm} \min_{X\in \mathcal{X}}\sum_{i=1}^4 {\rm vol}(D_i(C,X))\\
&\hspace{-0.2cm}\ge &\hspace{-0.2cm}\sum_{i=1}^4\min_{X\in \mathcal{X}}{\rm vol}(D_i(C,X))\\
&\hspace{-0.2cm}=&\hspace{-0.2cm} 4\cdot \mu_1 (C). \hspace{7.2cm}(3.5)
\end{eqnarray*}

Thus, we have proved the following result.

\medskip
\noindent
{\bf Lemma 3.2.}
$$\mu (C, {\bf e}_1)\ge 4\cdot \mu_1 (C).$$

\bigskip
Based on Lemma 3.1, we assume that $X$ is a discrete finite set of smallest cardinality such that $C+X$ is a packing in $\mathbb{E}^3$ and
$$\mu_1 (C)={\rm vol}(D_1(C,X)).\eqno (3.6)$$
By computing the value of ${\rm vol}(D_1(C,\Lambda_2))$, where $\Lambda_2$ is defined by (2.11), we obtain
$${\rm vol}(D_1(C,X))\le {\rm vol}(D_1(C,\Lambda_2))={1\over 6}\times {1\over 2}.\eqno (3.7)$$
For convenience, we write $X'=X\setminus \{ {\bf o}\}$ and enumerate the points of $X'$ with non-decreasing $x$-coordinates. Then ${\bf x}_1$ has the smallest $x$-coordinate among the points in $X'$ and
$$X'\subset \left\{ (x, y, z): 0<x<3,\ -1<y<2,\ -1<z<2\right\}\setminus {\rm int}(2C).\eqno (3.8)$$

If ${\bf x}_1=(x_1, y_1, z_1)\not\in \partial (2C)$ and $\epsilon $ is a positive number, we write ${\bf x}'_1 =(x_1-\epsilon , y_1, z_1)$ and
$$X_1={\bf x}'_1 \cup (X\setminus \{ {\bf x}_1\}).$$
When $\epsilon $ is small one can deduce that $C+X_1$ is a packing in $\mathbb{E}^3$ and
$${\rm vol}(D_1(C,X_1))<{\rm vol}(D_1(C,X )),$$
which contradicts the minimal assumption on ${\rm vol}(D_1(C,X))$. Thus we get
$${\bf x}_1\in \partial (2C).$$
In other words, $C+{\bf x}_1$ touches $C$ at its boundary.

For convenience, we write $x_0={\bf o}$ and define
\begin{eqnarray*}
F_0&\hspace{-0.2cm}=&\hspace{-0.2cm}\{ (x,y,z):\ x=1,\ |y|+|z|\le 1\},\\
F_1&\hspace{-0.2cm}=&\hspace{-0.2cm}\{ (x,y,z):\ x\ge 0,\ y\ge 0,\ z\ge 0,\ x+y+z=2\},\\
F_2&\hspace{-0.2cm}=&\hspace{-0.2cm}\{ (x,y,z):\ x\ge 0,\ y\le 0,\ z\ge 0,\ x-y+z=2\},\\
F_3&\hspace{-0.2cm}=&\hspace{-0.2cm}\{ (x,y,z):\ x\ge 0,\ y\le 0,\ z\le 0,\ x-y-z=2\},\\
F_4&\hspace{-0.2cm}=&\hspace{-0.2cm}\{ (x,y,z):\ x\ge 0,\ y\ge 0,\ z\le 0,\ x+y-z=2\}.
\end{eqnarray*}
In fact, $F_0$ is a square facet of $C$, and $F_1,$ $F_2$, $F_3$ and $F_4$ are four triangular facets surrounding $F_0$. Recall that if ${\bf x}_0{\bf x}_i\prec F_j$ holds for some ${\bf x}_i\in X'$, then the hyperplane generated by $F_j$ separates ${\rm int}(C)$ and ${\rm int}(C)+{\bf x}_i$. Thus, by the minimal assumption on the cardinality of $X$, one can deduced that
$${\bf x}_0{\bf x}_i\prec F_0\cup {\rm int}(F_1)\cup {\rm int}(F_2) \cup {\rm int}(F_4) \eqno (3.9)$$
holds for all points ${\bf x}_i\in X'$. In particular, $X'$ has a point ${\bf x}'=(x', y', z')$ satisfying ${\bf x}_0{\bf x}'\prec F_1$, $0<x'<2,$ $0<y'<2$ and $0<z'<2$. Otherwise, the tetrahedron with vertices $(1,1,1)$, $(0,1,1)$, $(1,0,1)$ and $(1,1,0)$ would be a subset of $D_1(C,X)$ and consequently
$${\rm vol}(D_1(C,X))\ge {1\over 6},$$
which contradicts to (3.7). Without loss of generality, we assume that ${\bf x}'$ has the smallest $x$-coordinate
among all points of this kind in $X$.

Let ${\bf x}_i=(x_i, y_i, z_i)\in X'$ be a point satisfying ${\bf x}_0{\bf x}_i\prec F_4$ (or ${\bf x}_0{\bf x}_i\prec F_2$), and let $D$ denote the regular hexagon with vertices $(1,0,1)$, $(1,1,0)$, $(0,1,-1)$, $(-1,0,-1)$, $(-1,-1,0)$ and $(0,-1,1)$. If $x'+y'+z'>4$, then $C+{\bf x}_i$ can not prevent $C+{\bf x}'$ from touching $C$ along the $-{\bf e}_1$ direction. In fact, if $C+{\bf x}_i$ blocks $C+{\bf x}'$ from moving in the $-{\bf e}_1$ direction at $F_1+{\bf x}_i$, then by (3.8) we get
$$\left\{ \begin{array}{ll}
x_i+y_i-z_i\ge 4&\\
(x'-x_i)+(y'-y_i)+(z'-z_i)=4&
\end{array}\right.$$
and therefore
$$x'+y'=4+x_i+y_i+(z_i-z')\ge 8+z_i+(z_i-z')\ge 4,$$
which contradicts the assumption on $y'<2$ and $z'<2$. If $C+{\bf x}_i$ blocks $C+{\bf x}'$ from moving in the $-{\bf e}_1$ direction by other part of $\partial (C)+{\bf x}_i$, then the hyperplane containing $D+{\bf x}_i$ will separate $D$ and $D+{\bf x}'$, which contradicts the assumption on ${\bf x}_0{\bf x}'\prec F_1$ and ${\bf x}_0{\bf x}_i\prec F_4$. Thus, defining ${\bf x}^*=(x'-\epsilon , y', z')$ and taking
$$X_2={\bf x}^* \cup (X\setminus \{ {\bf x}'\}),$$
when $\epsilon $ is a suitable small positive number one can deduce that $C+X_2$ is a packing in $\mathbb{E}^3$ and
$${\rm vol}(D_1(C,X_2))<{\rm vol}(D_1(C,X )),$$
which contradicts the minimal assumption on ${\rm vol}(D_1(C,X))$. As a conclusion we have proved the following assertion.

\medskip\noindent
{\bf Lemma 3.3.} {\it If $C+X$ is a packing such that $\mu_1(C)={\rm vol}(D_1(C,X))$, then there is a point
$(x, y, z)\in X$ satisfying $0<x<2$, $0<y<2$, $0<z<2$ and $x+y+z=4$. In other words, there is a translate $C+{\bf x}$ in $C+X$ which touches $C$ at some interior point of $F_1$.}

\medskip
Now we prove the following result, which is a preliminary estimate for $\mu_1(C)$.

\medskip
\noindent
{\bf Theorem 3.1.}
$$\mu_1(C)\ge {1\over 6}\times {{4\sqrt{2}+2}\over {25+22\sqrt{2}}}.$$

\medskip\noindent
{\bf Proof.} Clearly $F_1$ is the triangular facet of $C$ with vertices $(0,1,1)$, $(1,0,1)$ and $(1,1,0)$.
Let $T_0$ denote the orthogonal tetrahedron with vertices $(0,0,0)$, $(1,0,0)$, $(0,1,0)$ and
$(0,0,1)$, and let $X$ be a discrete set such that
$$\mu_1(C)={\rm vol}(D_1(C,X)).$$

It follows by Lemma 3.3 that there is a point ${\bf x}_1=(x_1, y_1, z_1)\in X$ such that $C+{\bf x}_1$ touches $C$ at the interior of $F_1$. Then $C\cap (C+{\bf x}_1)$ is either a centrally symmetric hexagon or a parallelogram.

\begin{figure}[hc]
\includegraphics[height=5cm,width=7.5cm,angle=0]{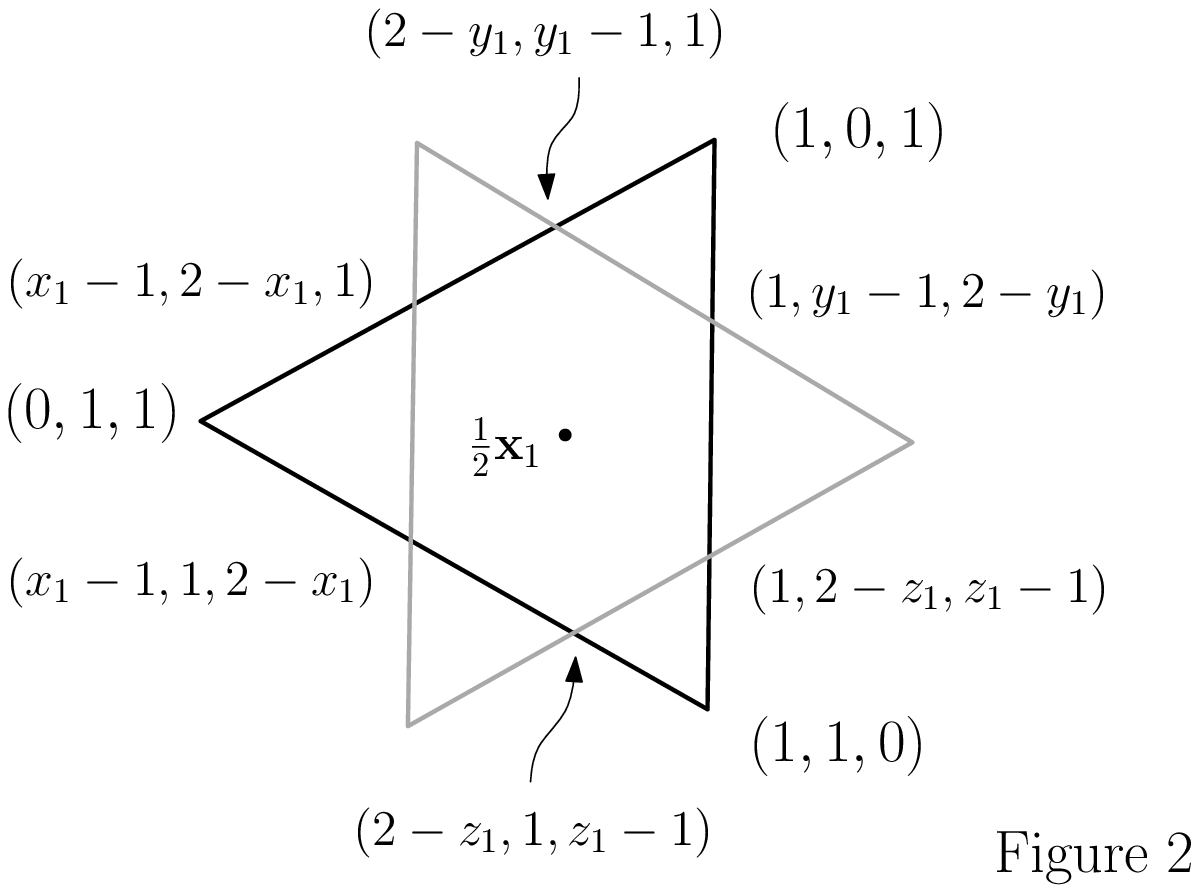}
\end{figure}

When $C\cap (C+{\bf x}_1)$ is an hexagon, its vertices are $(1, 2-z_1, z_1-1)$, $(1, y_1-1, 2-y_1),$ $(2-y_1, y_1-1, 1)$, $(x_1-1, 2-x_1, 1)$, $(x_1-1, 1, 2-x_1)$ and $(2-z_1, 1, z_1-1)$, as shown
in Figure 2, where
$$x_1+y_1+z_1=4.\eqno (3.10)$$
Let $T_1$ denote the tetrahedron with vertices $(x_1-1,1,1)$, $(0,1,1)$, $(x_1-1,2-x_1,1)$ and $(x_1-1,1,2-x_1)$,
let $T_2$ denote the tetrahedron with vertices $(1,y_1-1,1)$, $(2-y_1,y_1-1,1)$, $(1,0,1)$ and $(1,y_1-1,2-y_1)$,
let $T_3$ denote the tetrahedron with vertices $(1,1,z_1-1)$, $(2-z_1,1,z_1-1)$, $(1,2-z_1,z_1-1)$ and $(1,1,0)$,
and let $T_4$ denote the tetrahedron with vertices $(1,y_1-1,z_1-1)$, $(x_1,y_1-1,z_1-1)$, $(1,2-z_1,z_1-1)$ and $(1,y_1-1,2-y_1)$. Clearly, all $T_1$, $T_2$, $T_3$ and $T_4$ are homothetic to $T_0$ with ratios $x_1-1$, $y_1-1$, $z_1-1$ and $x_1-1$, respectively, and
$${\rm int}(T_i)\cap {\rm int}(T_j)=\emptyset $$
holds whenever $i\not= j$. By routine arguments based on Lemmas 2.3 and 2.4 we get
$$T_i\subseteq D_1(C,X),\qquad i=1, 2, 3, 4.$$
Thus, by (3.10), we get
\begin{eqnarray*}
\hspace{4.2cm}{\rm vol}(D_1(C,X))&\hspace{-0.2cm}\ge &\hspace{-0.2cm}{\rm vol}(T_1\cup T_2\cup T_3\cup T_4)\\
&\hspace{-0.2cm}=&\hspace{-0.2cm}{1\over 6} \left(2(x_1-1)^3+(y_1-1)^3+(z_1-1)^3\right)\\
&\hspace{-0.2cm}\ge &\hspace{-0.2cm} {1\over 6} \left(2(x_1-1)^3+2\left(1-{1\over 2}x_1\right)^3\right)\\
&\hspace{-0.2cm}\ge &\hspace{-0.2cm}{1\over 6}\times {{4\sqrt{2}+2}\over {25+22\sqrt{2}}}, \hspace{6.1cm}(3.11)
\end{eqnarray*}
where the last equality holds if and only if $x_1=1+{1\over {2\sqrt{2}+1}}$.

When $C\cap (C+{\bf x}_1)$ is a parallelogram, by similar arguments, it can be proved that
$${\rm vol}(D_1(C,X))\ge {1\over 6}\times {1\over 4}.\eqno (3.12)$$

As a conclusion of (3.11) and (3.12), noticing that
$${{4\sqrt{2}+2}\over {25+22\sqrt{2}}}\approx 0.1364549\cdots ,$$
Theorem 3.1 is proved. \hfill{$\Box$}

\medskip
By Lemma 2.1, Lemma 3.2, Theorem 3.1 and (2.2) it follows that

\medskip
\noindent
{\bf Corollary 3.1.}
$$\delta^t(C)\le 0.9865382\cdots \qquad and \qquad \delta^t(T)\le 0.3946153\cdots .$$

\vspace{0.9cm}\noindent
{\large\bf 4. Cuboctahedral Packings, A Detailed Computation}

\vspace{0.4cm}\noindent
In this section we prove the following result.

\medskip\noindent
{\bf Theorem 4.1.}
$$\mu_1 (C)={1\over 6}\left({5\over 9}-{4\over 9}\sqrt{1\over {10}}\right).$$

\bigskip
\noindent
{\bf Remark 4.1.} {\it We notice that} $${5\over 9}-{4\over 9}\sqrt{1\over {10}}\approx 0.415009881\cdots . $$

\medskip
For convenience, we write
$$\alpha ={1\over 3}\left(2+4 \cos {{\pi +\arctan \sqrt{63} }\over 3}\right)\approx 0.7223517\cdots $$
and
$$\beta =1+\sqrt[3]{{1\over 2}}\approx 1.7937005\cdots .$$
In fact, they are roots of the equations $3x^3-6x^2+2=0$ and $(x-1)^3-{1\over 2}=0$, respectively.

\begin{figure}[hc]
\includegraphics[height=4.2cm,width=5.5cm,angle=0]{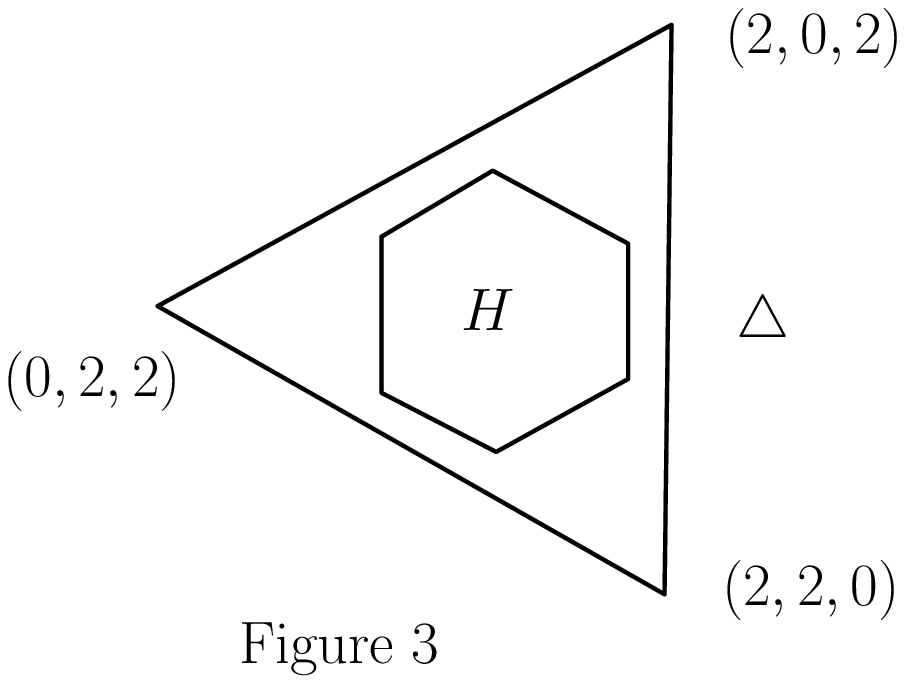}
\end{figure}

Let $\triangle $ denote the triangle with vertices $(0,2,2)$, $(2,0,2)$ and $(2,2,0)$, as illustrated in Figure 3, and write
$$H=\{ (x,y,z):\ \alpha \le x\le \beta,\ \alpha \le y\le \beta,\ \alpha \le z\le \beta,\ x+y+z=4\}.$$

\medskip\noindent
{\bf Lemma 4.1.} {\it If ${\bf x}_1=(x_1, y_1, z_1)\in \triangle\setminus H$, we have}
$${\rm vol}(D_1(C, X))\ge {1\over 6}\times {1\over 2}.$$

\medskip\noindent
{\bf Proof.} When $x_1\ge \beta $, by Lemma 2.3 and Lemma 2.4 it can be shown that the tetrahedron $T_1$ with vertices $(x_1-1, 1,1)$, $(0,1,1)$, $(x_1-1, 2-x_1,1)$ and $(x_1-1, 1, 2-x_1)$ is a subset of $D_1(C,X)$. The tetrahedron $T_1$ is homothetic to the orthogonal tetrahedron $T_0$ defined in the proof of Theorem 3.1 with ratio $x_1-1$. Thus, we get
$${\rm vol}(D_1(C,X))\ge {\rm vol}(T_1)={1\over 6} (x_1-1)^3\ge {1\over 6}\left(\sqrt[3]{1\over 2}\right)^3={1\over 6}\times {1\over 2}.\eqno (4.1)$$
The cases $y_1\ge \beta $ and $z_1\ge \beta $ can be treated by similar arguments.

When $x_1\le \alpha $, by Lemmas 2.3 and 2.4 it can be verified that both tetrahedra $T_2$ and $T_3$ are contained in $D_1(C,X)$, where $T_2$ has vertices $(1,y_1-1,1)$, $(2-y_1, y_1-1, 1)$, $(1,0,1)$ and $(1,y_1-1,2-y_1)$, $T_3$ has vertices $(1,1,z_1-1)$, $(2-z_1,1,z_1-1)$, $(1,2-z_1,z_1-1)$ and $(1,1,0)$. Clearly all $T_2$, $T_3$ and $T_2\cap T_3$ are homothetic to $T_0$ with ratios $y_1-1$, $z_1-1$ and $1-x_1$, respectively. Thus, recalling the definition of $\alpha $, we have
\begin{eqnarray*}
\hspace{3.3cm}{\rm vol}(D_1(C,X))&\hspace{-0.2cm}\ge &\hspace{-0.2cm}{\rm vol}(T_2\cup T_3)={\rm vol}(T_2)+{\rm vol}(T_3)-{\rm vol}(T_2\cap T_3)\\
&\hspace{-0.2cm}=&\hspace{-0.2cm}{1\over 6}\left((y_1-1)^3+(z_1-1)^3-(1-x_1)^3\right)\\
&\hspace{-0.2cm}\ge &\hspace{-0.2cm} {1\over 6} \left(2\left(1-{1\over 2}x_1\right)^3-(1-x_1)^3\right)\\
&\hspace{-0.2cm}\ge &\hspace{-0.2cm} {1\over 6} \left(2\left(1-{1\over 2}\alpha \right)^3-(1-\alpha )^3\right)\\
&\hspace{-0.2cm} = &\hspace{-0.2cm} {1\over 6}\times {1\over 2}. \hspace{8.7cm}(4.2)
\end{eqnarray*}
The cases $y_1\le \alpha $ and $z_1\le \alpha $ can be proved by similar arguments.

As a conclusion of (4.1) and (4.2) the lemma is proved. \hfill{$\Box $}

\medskip
Assume that $X$ is a discrete set of points such that $C+X$ is a packing in $\mathbb{E}^3$ satisfying
$${\rm vol}(D_1(C,X))=\mu_1(C).$$
Let $\triangle $ denote the triangle with vertices $(0,2,2)$, $(2,0,2)$ and $(2,2,0)$, let $\triangle^*$ denote the triangle with vertices $(2,1,1)$, $(1,2,1)$ and $(1,1,2)$, let $\triangle_1$ denote the triangle with vertices $(2,0,2)$, $(2,0.5, 1.5)$ and $(1,1,2)$, let $\triangle_2$ denote the triangle with vertices $(2,0.5,1.5)$, $(2,1,1)$ and $(1,1,2)$, let $\triangle_3$ denote the triangle with vertices $(2,1,1)$, $(1,1,2)$ and $(1,1.5,1.5)$, let
$\triangle_4$ denote the triangle with vertices $(1, 1.5, 1.5)$, $(1,1,2)$ and $(0,2,2)$, and let $\triangle '$ denote the triangle with vertices $(0,1,1)$, $(1,0,1)$ and $(1,1,0)$.

\begin{figure}[hc]
\includegraphics[height=5cm,width=6.3cm,angle=0]{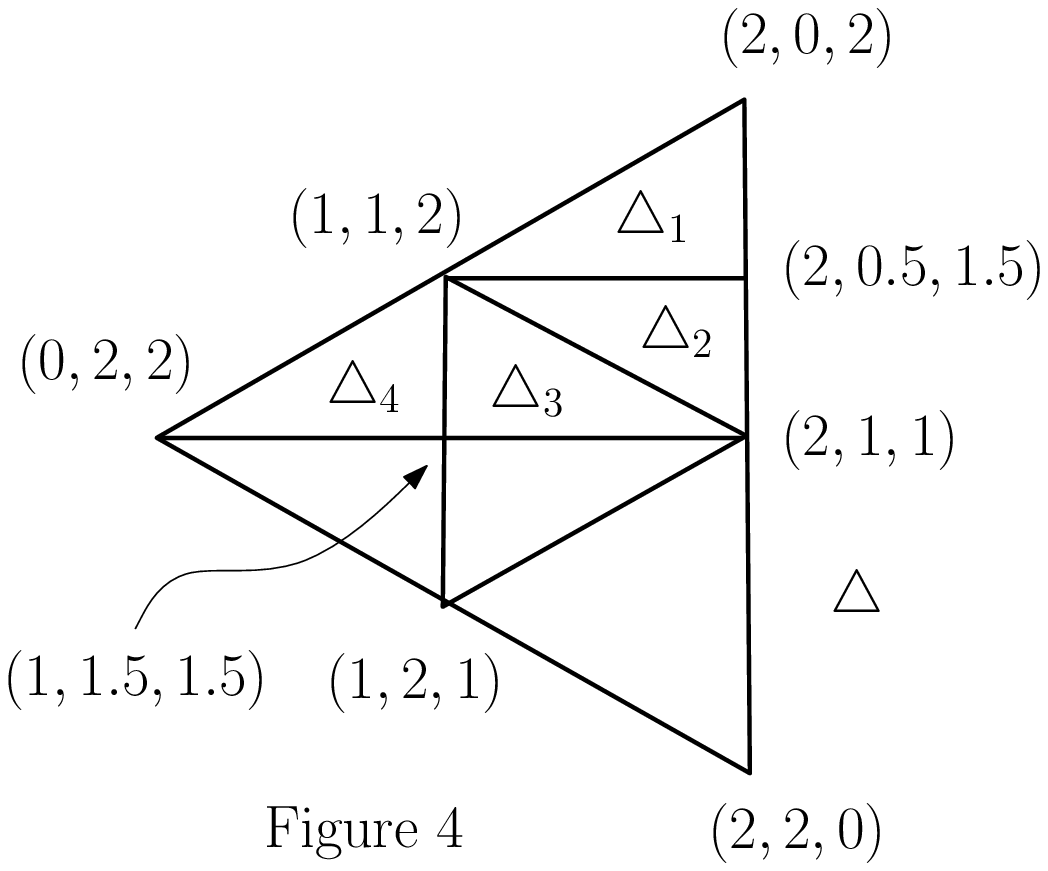}
\end{figure}

To estimate $\mu_1(C)$, based on Lemma 3.3, we assume that
$${\bf x}_1=(x_1, y_1, z_1)\in {\rm int}(\triangle )\cap X$$
and therefore
$$x_1+y_1+z_1=4.\eqno (4.3)$$
It can be shown that $C\cap (C+{\bf x}_1)$ is a centrally symmetric hexagon if ${\bf x}_1\in \triangle^*$ and is a parallelogram if ${\bf x}_1\in {\rm int}(\triangle )\setminus \triangle^*$. By symmetry, we consider four cases with respect to ${\bf x}_1\in \triangle_i$ by Lemmas 4.2-4.5, respectively.

\medskip\noindent
{\bf Lemma 4.2.} {\it If ${\bf x}_1\in \triangle_1$, we have}
$${\rm vol}(D_1(C,X))\ge {1\over 6}\times {9\over {16}}.$$

\medskip\noindent
{\bf Proof.} In this case $C+{\bf x}_1$ touches $C$ at a parallelogram with vertices $(1,0,1)$, $(1, 2-z_1, z_1-1)$, $(x_1-1, y_1, z_1-1)$ and $(x_1-1,2-x_1,1)$. Then, by Lemma 2.3 and Lemma 2.4 it can be shown that $D_1(C,X)$ contains three tetrahedra $T_1$, $T_3$ and $T_5$ and a prismoid $P$, where $T_1$ has vertices $(x_1-1, 1,1)$, $(0,1,1)$, $(x_1-1, 2-x_1,1)$ and $(x_1-1, 1, 2-x_1)$, $T_3$ has vertices $(1,1, z_1-1)$, $(2-z_1, 1, z_1-1)$, $(1, 2-z_1, z_1-1)$ and $(1,1,0)$, $T_5$ has vertices $(1,0, z_1-1)$, $(3-z_1, 0, z_1-1)$, $(1, 2-z_1, z_1-1)$ and $(1,0,1)$, and $P$ has six vertices $(x_1-1, y_1, z_1-1)$, $(x_1-1, 1, z_1-1)$, $(x_1-y_1, 1, z_1-1)$, $(x_1-1, 1, 1)$, $(x_1-1, 2+y_1-z_1, 1)$ and $(2-2y_1, 1,1)$. In other words,
$$P=\{(x,y,z):\ x\ge x_1-1,\ y\le 1,\ z_1-1\le z\le 1,\ (x-x_1)-(y-y_1)+(z-z_1)\le -2 \}.$$

Clearly, all $T_1$, $T_3$, $T_5$ and $T_1\cap T_3$ are homothetic to the orthogonal tetrahedron $T_0$ which has vertices $(0,0,0)$, $(1,0,0)$, $(0,1,0)$ and $(0,0,1)$ with ratios $x_1-1$, $z_1-1$, $2-z_1$ and $1-y_1$, respectively. The prismoid $P$ can be obtained by cutting off a tetrahedron $T_7$ with vertices $(x_1-1, 1,1)$, $(2-2y_1, 1,1)$, $(x_1-1, 2+y_1-z_1, 1)$ and $(x_1-1, 1, z_1-y_1)$ from a tetrahedron $T_6$ with vertices $(x_1-1, 1, z_1-1)$, $(x_1-y_1, 1, z_1-1)$, $(x_1-1, y_1, z_1-1)$ and $(x_1-1, 1, z_1-y_1)$. Both $T_6$ and $T_7$
are homothetic to $T_0$ with ratios $1-y_1$ and $z_1-y_1-1$, respectively.

By a routine computation based on the assumption $(x_1, y_1, z_1)\in \triangle_1$ one can deduce
\begin{eqnarray*}
{\rm vol}(T_1\cup T_3\cup T_5\cup P)&\hspace{-0.2cm}=&\hspace{-0.2cm}{1\over 6} \left((x_1-1)^3+(z_1-1)^3+(2-z_1)^3-(1-y_1)^3 +(1-y_1)^3-(z_1-y_1-1)^3\right)\\
&\hspace{-0.2cm}=&\hspace{-0.2cm}{1\over 6} \left((x_1-1)^3+(z_1-1)^3+(2-z_1)^3-(2z_1+x_1-5)^3\right)\\
&\hspace{-0.2cm}\ge &\hspace{-0.2cm}{1\over 6} \left(\left({1\over 2}\right)^3+\left({3\over 4}\right)^3+\left({1\over 4}\right)^3\right)\\
&\hspace{-0.2cm}= &\hspace{-0.2cm} {1\over 6}\times {9\over {16}},
\end{eqnarray*}
where the equality holds if and only if ${\bf x}_1=({3\over 2}, {3\over 4}, {7\over 4}).$ The inequality can be deduced by the fact that the function $f(x_1, z_1)=(x_1-1)^3+(z_1-1)^3+(2-z_1)^3-(2z_1+x_1-5)^3$ attains its minimum on the boundary of $\triangle_1$. Thus, in this case we have
$${\rm vol}(D_1(C,X))\ge {\rm vol}(T_1\cup T_3\cup T_5\cup P)\ge {1\over 6}\times {9\over {16}}.$$
The lemma is proved. \hfill{$\Box $}

\medskip\noindent
{\bf Remark 4.2.} {\it It can be verified that ${\bf x}_1=({3\over 2}, {3\over 4}, {7\over 4})\in {\rm int}(H)$, where $H$ was defined above Lemma $4.1$. Therefore Lemma $4.2$ can not be covered by Lemma $4.1$.}

\medskip
\noindent
{\bf Lemma 4.3.} {\it If ${\bf x}_1\in \triangle_2$, we have}
$${\rm vol}(D_1(C,X))\ge {1\over 6}\left(6-12\sqrt{6\over {29}}\right).$$

\medskip\noindent
{\bf Proof.} In this case $C+{\bf x}_1$ touches $C$ at the parallelogram with vertices
$(1,0,1)$, $(1, 2-z_1, z_1-1)$, $(x_1-1, y_1, z_1-1)$ and $(x_1-1,2-x_1,1)$. Then, by Lemmas 2.3 and 2.4 it can be shown that  $D_1(C,X)$ contains all $T_1$, $T_3$, $T_5$ and $T_6$, where $T_1$ is the tetrahedron with vertices $(x_1-1,1,1)$, $(0,1,1)$, $(x_1-1,2-x_1,1)$ and $(x_1-1,1,2-x_1)$, $T_3$ is the tetrahedron with vertices $(1,1,z_1-1)$, $(2-z_1,1,z_1-1)$, $(1,2-z_1,z_1-1)$ and $(1,1,0)$, $T_5$ is the tetrahedron with vertices $(1,0, z_1-1)$, $(3-z_1, 0, z_1-1)$, $(1, 2-z_1, z_1-1)$ and $(1,0,1)$, and $T_6$ is the tetrahedron with vertices $(x_1-1, 1, z_1-1)$, $(x_1-y_1, 1, z_1-1)$, $(x_1-1, y_1, z_1-1)$ and $(x_1-1, 1, z_1-y_1)$.

Clearly, all $T_1$, $T_3$, $T_5$, $T_6$ and $T_1\cap T_3$ are homothetic to the orthogonal tetrahedron $T_0$ with ratios $x_1-1$, $z_1-1$, $2-z_1$, $1-y_1$ and $1-y_1$, respectively. Thus, we have

\begin{eqnarray*}
\hspace{1.2cm}{\rm vol}(T_1\cup T_3\cup T_5\cup T_6)&\hspace{-0.2cm}=&\hspace{-0.2cm}
{1\over 6} \left((x_1-1)^3+(z_1-1)^3-(1-y_1)^3+(2-z_1)^3+(1-y_1)^3\right)\\
&\hspace{-0.2cm}=&\hspace{-0.2cm}{1\over 6} \left((x_1-1)^3+(z_1-1)^3+(2-z_1)^3\right).
\hspace{4.9cm} (4.4)\end{eqnarray*}

Let $U$ denote the rectangle with vertices $(1,0,0)$, $(1,0,z_1-1)$, $(1,1,z_1-1)$ and $(1,1,0)$ and define
$$G_1(C,X)=\Bigl(\bigcup_{{\bf x}\in U}s(C,X,{\bf e}_1, {\bf x})\Bigr)\cap 2W.\eqno (4.5)$$
We proceed to estimate the volume of $G_1(C,X)$. To this end, let $S_0$ to be the halfspace $\{ (x, y, z):\ x\ge 1\}$, let $S_1$ to be the halfspace $\{ (x, y, z):\ y\ge 0\}$, let $S_2$ to be the halfspace $\{ (x, y, z):\ z\ge 0\}$, let $S_3$ to be the halfspace $\{ (x, y, z):\ x\le 2\}$, let $S_4$ to be the halfspace $\{ (x, y, z):\ y\le 1\}$, let $S_5$ to be the halfspace $\{ (x, y, z):\ z\le z_1-1\}$, let $S_6$ to be the halfspace $\{ (x, y, z):\ x+y-z\le 2\}$, let $S_7$ to be the halfspace $\{ (x, y, z):\ (x-x_1)-(y-y_1)-(z-z_1)\le 2\}$, and define
$$S=\bigcap_{i=0}^7S_i.$$
Furthermore, we write ${\bf x}_t =(2, t , z_1-2)$ and ${\bf x}'_t =(2, t -2 , z_1-2)$ and define
$$P(t) =S\setminus \left\{ (C+{\bf x}_t )\cup (C+{\bf x}'_t )\right\}.$$

For convenience, we write
$$F_5=\{ (x,y,z):\ z=-1,\ |x|+|y|\le 1\}.$$
If ${\bf x}\in X$ and ${\rm int}(S)\cap (C+{\bf x})\not= \emptyset $, by reducing the $x$-coordinate of ${\bf x}$ until $C+{\bf x}$ touches $C$ one can deduce ${\bf x}_0{\bf x}\prec F_0$. In addition, if there are more than one ${\bf x}\in X$ satisfying ${\rm int}(S)\cap (C+{\bf x})\not=\emptyset $, they can not block each others in this reducing process. To see this, if $C+{\bf x}_2$ touches $C$ at some relative interior points of $F_0$, and $C+{\bf x}_2$ blocks $C+{\bf x}_3$ from moving in $-{\bf e}_1$ direction by $F_i+{\bf x}_2$, one can deduce contradictions by considering five cases with respect to $i=0, 1, 2, 3$ and $4$.

The $i=0$ case is obvious. When $i=1$ or $2$, we have ${\bf x}_1{\bf x}_3\prec F_5+{\bf x}_1$, $z_3\le z_1-2$, $z_2+1\le z_3,$ and therefore
$$z_2\le z_3-1\le z_1-3\le -1,$$
which contradicts the assumption that ${\rm int}(S)\cap (C+{\bf x}_2)\not=\emptyset .$ When $i=3$, we get
$$\left\{ \begin{array}{ll}
|y_2|+|z_2|\le 2 &\\
(x_3-2)-(y_3-y_2)-(z_3-z_2)=4&
\end{array}\right.$$
and therefore
$$x_3-y_3-z_3=6-(y_2+z_2)\ge 4,$$
which contradicts the assumption that ${\rm int}(S)\cap (C+{\bf x}_3)\not=\emptyset .$ When $i=4$, we get
$$\left\{ \begin{array}{ll}
|y_2|+|z_2|\le 2 &\\
(x_3-2)+(y_3-y_2)-(z_3-z_2)=4&
\end{array}\right.$$
and therefore
$$x_3+y_3-z_3=6+y_2-z_2\ge 4,$$
which contradicts the assumption that ${\rm int}(S)\cap (C+{\bf x}_3)\not=\emptyset .$

Then, by increasing the $z$-coordinate of ${\bf x}$ until $C+{\bf x}$ touches $C+{\bf x}_1$, it can be shown that
$${\rm vol}(G_1(C,X))\ge \min_t\{ {\rm vol}(P(t))\}.\eqno (4.6)$$
\begin{figure}[hc]
\includegraphics[height=7.8cm,width=8.6cm,angle=0]{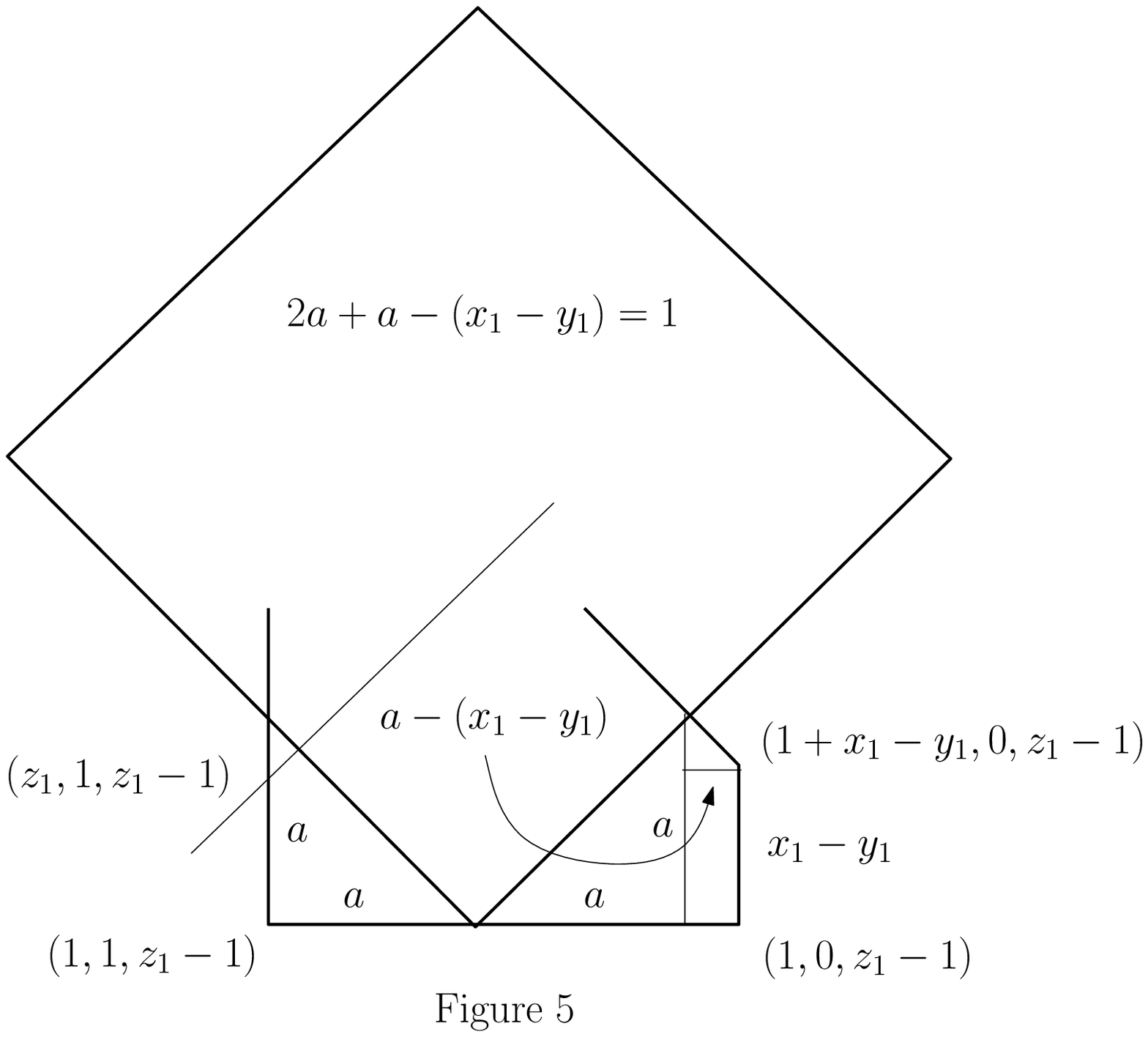}
\end{figure}

It can be shown that the three hyperplanes $y=0$, $z=z_1-1$ and $(x-x_1)-(y-y_1)-(z-z_1)=2$ intersect at
$(1+x_1-y_1, 0, z_1-1)$, and the three hyperplanes $y=1$, $z=z_1-1$ and $x+y-z=2$ meet at $(z_1, 1, z_1-1)$.
Based on Figure 5, which illustrates a possible $C+{\bf x}$ touching $C+{\bf x}_1$ in the half plane $\{ (x,y,z):\ z=z_1-1,\ x\ge 1\}$, by studying the derivative of ${\rm vol}(P(t))$ in terms of parts of surface areas of $P(t)$ it can be shown that
$${\rm vol}(P(t))\ge {1\over 6} \left\{\begin{array}{ll}
2\left({1\over 2}\right)^3& \mbox{ if $x_1-y_1\ge {1\over 2}$, $z_1-1\ge {1\over 2}$,}\\
2\left(z_1-1\right)^3& \mbox{ if $x_1-y_1\ge {1\over 2}$, $z_1-1\le {1\over 2}$,}\\
2\left({{1+x_1-y_1}\over 3}\right)^3 & \mbox{ if $x_1-y_1\le {1\over 2}$, $z_1-1\ge {1\over 3}(1+x_1-y_1)$,}\\
2\left(z_1-1\right)^3 & \mbox{ if $x_1-y_1\le {1\over 2}$, $z_1-1\le {1\over 3}(1+x_1-y_1)$.}
\end{array}\right.$$
Thus, we get
$$\min_{t}\{{\rm vol}(P(t))\}\ge {1\over 6}\min \left\{ 2\left( {{1+x_1-y_1}\over 3}\right)^3,\
2 (z_1-1)^3,\ 2\left({1\over 2}\right)^3 \right\}.\eqno (4.7)$$

By routine computations it can be shown that, when ${\bf x}_1\in \triangle_2$,
$$(x_1-1)^3+(z_1-1)^3+(2-z_1)^3+ 2\left( {{1+x_1-y_1}\over 3}\right)^3\ge 6-12\sqrt{{6\over
{29}}},\hspace{0.2cm}\eqno (4.8)$$
$$(x_1-1)^3+(z_1-1)^3+(2-z_1)^3+ 2(z_1-1)^3\ge {6\over {5+2\sqrt{6}}}\ge 6-12\sqrt{{6\over
{29}}}\eqno (4.9)$$
and
$$(x_1-1)^3+(z_1-1)^3+(2-z_1)^3+ 2\left( {1\over 2}\right)^3\ge {{11+2\sqrt{2}}\over {12+8\sqrt{2}}}
\ge 6-12\sqrt{{6\over
{29}}}.\eqno (4.10)$$
Thus, in this case it can be deduced from (4.4), (4.6)-(4.10) that
\begin{eqnarray*}
{\rm vol}(D_1(C, X)) &\hspace{-0.2cm}\ge &\hspace{-0.2cm} {\rm vol}(T_1\cup T_3\cup T_5\cup T_6)+{\rm vol}(G_1(C,X))\\
&\hspace{-0.2cm}\ge &\hspace{-0.2cm} {\rm vol}(T_1\cup T_3\cup T_5\cup T_6)+\min_t\{ {\rm vol}(P(t))\}\\
&\hspace{-0.2cm}\ge &\hspace{-0.2cm} {1\over 6} \left((x_1-1)^3+(z_1-1)^3+(2-z_1)^3+ \min \left\{ 2\left( {{1+x_1-y_1}\over 3}\right)^3,\
2 (z_1-1)^3,\ 2\left({1\over 2}\right)^3 \right\}\right)\\
&\hspace{-0.2cm}\ge &\hspace{-0.2cm}{1\over 6}\left(6-12\sqrt{{6\over {29}}}\right).\end{eqnarray*}
Lemma 4.3 is proved. \hfill{$\Box $}

\bigskip\noindent
{\bf Remark 4.3.} {\it We notice that} $$6-12\sqrt{{6\over {29}}}\approx 0.541694086\cdots .$$

\bigskip\medskip
We define
\begin{eqnarray*}
P_1\hspace{-0.25cm}&=\hspace{-0.25cm}&\{ (x,y,z):\ 0\le x< 2,\ 0< y<1,\ 0<z<1,\ x-y+z< 2,\
x-y-z<1,\ x+y-z<2\},\\
P_2\hspace{-0.25cm}&=\hspace{-0.25cm}&\{ (x,y,z):\ 1\le x < 3,\ -1<y<2,\ -1<z<2,\ x-y+z<4,\ x-y-z<3,\ x+y-z<4\}
\end{eqnarray*}
and
$$D'_1(C,X)=D_1(C,X)\cap P_1.$$

It follows by Lemma 3.3 that $X$ has a point ${\bf x}_1=(x_1, y_1, z_1)$ which belongs to ${\rm int}(\triangle )$.
For convenience, we write ${\bf x}_0={\bf o}$ and define
$$X'=(X\cap P_2)\cup \{ {\bf x}_0, {\bf x}_1\}.\eqno (4.11)$$
It can be verified by Lemmas 2.3 and 2.4 that
$$D'_1(C,X')=D'_1(C, X).$$

In the forthcoming proofs of both Lemma 4.4 and Lemma 4.5 we will estimate the minimal value of ${\rm vol}(D'_1(C,X'))$ instead of ${\rm vol}(D_1(C,X))$. Thus we may assume that
$$\biggl( (C+{\bf x}_1)\cup C\biggr) \bigcap \biggl( \bigcup_{{\bf x}\in X'\setminus \{ {\bf x}_0,
{\bf x}_1\}}(C+{\bf x})\biggr)\not=\emptyset .\eqno (4.12)$$
For convenience, we write
$$J=\biggl( (C+{\bf x}_1)\cup C\biggr) \bigcap \biggl( \bigcup_{{\bf x}\in X'\setminus \{ {\bf x}_0, {\bf x}_1\}}(C+{\bf x})\biggr).$$

Let $F_i$ be the facets defined just above (3.9). In particular, $F_0$ is the square facet $\{ (x,y,z):\ x=1,\ |y|+ |z|\le 1\}$ of $C$, and $F_3$ is the triangular facet $\{ (x,y,z):\ \max \{ |x|,\ |y|,\ |z|\}\le 1,\ x-y-z=2\}$. It follows by (4.11) and (4.12) that
$$J\cap \left( {\rm int}(F_0)\cup \{ {\rm int}(F_3)+{\bf x}_1\} \cup \{ F_0+{\bf x}_1\}\right)\not= \emptyset.\eqno (4.13)$$

\medskip
\noindent
{\bf Lemma 4.4.} {\it If ${\bf x}_1\in \triangle_4\cap H$, then we have}
$${\rm vol}(D_1(C,X))\ge {1\over 6} \left({5\over 9}-{4\over 9}\sqrt{1\over {10}}\right).$$

\medskip\noindent
{\bf Proof.} In this case $C+{\bf x}_1$ touches $C$ at a parallelogram with vertices $(0,1,1)$, $(2-z_1, 1, z_1-1)$,
$(x_1, y_1-1, z_1-1)$ and $(2-y_1, y_1-1, 1)$. Then, by routine arguments based on Lemmas 2.3 and 2.4 it can be shown that $D'_1(C,X')$ contains both $T_2$ and $T_3$, where $T_2$ is the tetrahedron with vertices $(1, y_1-1, 1)$, $(2-y_1, y_1-1, 1)$, $(1,0,1)$ and $(1, y_1-1, 2-y_1)$, and $T_3$ is the tetrahedron with vertices $(1,1,z_1-1)$, $(2-z_1,1,z_1-1)$, $(1,2-z_1,z_1-1)$ and $(1,1,0)$. Clearly all $T_2$, $T_3$ and $T_2\cap T_3$ are homothetic to the orthogonal tetrahedron $T_0$ with ratios $y_1-1$, $z_1-1$ and $1-x_1$, respectively. Thus we have
$${\rm vol}(T_2\cup T_3)={1\over 6}\ ((y_1-1)^3+(z_1-1)^3-(1-x_1)^3)\ge {1\over 6}\times {1\over 4}.\eqno (4.14)$$

To estimate ${\rm vol}(D'_1(C,X'))$, based on (4.13) we consider three cases.

\medskip\noindent
{\bf Case 1.} $J\cap ({\rm int}(F_0)\cup ({\rm int}(F_3)+{\bf x}_1))=\emptyset .$ By the assumption on $X'$ we have
$$x\ge x_1+2$$
for all ${\bf x}=(x,y,z)\in X'\setminus \{ {\bf x}_0, {\bf x}_1\}$. Let $T^*$ denote the tetrahedron with
vertices $(1,0,0)$, $(1,1,0)$, $(1,0,1)$ and $(1+x_1, {1\over 2}, {1\over 2})$. Clearly we have
$$T^*\subset D'_1(C,X')$$
and
$${\rm vol}(T^*)={1\over 6}\ x_1.$$
Thus, applying (4.14) and (4.3), we get
\begin{eqnarray*}
\hspace{3.6cm}{\rm vol}(D'_1(C,X'))&\hspace{-0.2cm}\ge &\hspace{-0.2cm} {\rm vol}(T_2\cup T_3\cup T^*)\\
&\hspace{-0.2cm}= &\hspace{-0.2cm}{1\over 6}\ ((y_1-1)^3+(z_1-1)^3-(1-x_1)^3+x_1)\\
&\hspace{-0.2cm}\ge &\hspace{-0.2cm} {1\over 6}\left( 2\left(1-{{x_1}\over 2}\right)^3-(1-x_1)^3+x_1\right) \\
&\hspace{-0.2cm}\ge &\hspace{-0.2cm}{1\over 6}\left( 2\left( 1-{\alpha \over 2}\right)^3-(1-\alpha )^3+\alpha \right)\\
&\hspace{-0.2cm}> &\hspace{-0.2cm} {1\over 6}\times 1.\hspace{8.2cm} (4.15)
\end{eqnarray*}

\medskip\noindent
{\bf Case 2.} ${\rm int}(F_3+{\bf x}_1)\cap J\not=\emptyset $. Assume that $C+{\bf x}_2$, where ${\bf x}_2=(x_2,y_2,z_2)\in X'$, touches $C+{\bf x}_1$ at some interior points of $F_3+{\bf x}_1$. Clearly we have
$$y_2\ge y_1-2 \qquad {\rm and} \qquad z_2\ge z_1-2\ge y_1-2.$$
Then we consider two subcases.

\medskip\noindent
{\bf Subcase 2.1.} $y_2\ge z_1-2$. Let us define $T_8$ to be the tetrahedron with vertices $(1,1,z_1-1)$, $(2+x_1-y_1, 1, z_1-1)$, $(1, y_1-x_1, z_1-1)$ and $(1,1,2-2x_1)$, and define $T_9$ to be the
tetrahedron with vertices $(1, y_1-1, 1)$, $(2+x_1-z_1, y_1-1,1)$, $(1,2-2x_1,1)$ and $(1, y_1-1, z_1-x_1)$.

If there is a point ${\bf x}=(x,y,z)$ satisfying both
$${\rm int}(T_8)\cap (C+{\bf x} )\not= \emptyset $$
and
$$({\rm int}(C)+{\bf x})\cap (C+{\bf x}_i)=\emptyset , \quad i=0,1,2,$$
by Lemma 2.3 and Lemma 2.4 it can be verified that neither $C+{\bf x}_2$ nor $C+{\bf x}_1$ can block $C+{\bf x}$ from moving in $-{\bf e}_1$ direction. Thus, one can assume that ${\bf x}=(2,y,z)$ and therefore $C+{\bf x}$ touches $C$ at some interior points of $F_0$. Then the point $(1, y_2+1, z_1-2)$, which is one unit bellow $(1, y_2+1, z_1-1)$
in $z$ direction, and the origin ${\bf o}$ are on the same side of the hyperplane $\{ (x,y,z):\ x+y-z=2 \}$. Therefore we get $$y_2<2+(z_1-2)-1-1=z_1-2,$$ which contradicts the assumption on $y_2$. Thus we have
$${\rm int}(T_i)\cap (C+X')=\emptyset,\quad i=8,\ 9.$$
Clearly, both $T_8$ and $T_9$ are homothetic to $T_0$, with ratios $1+x_1-y_1$ and $1+x_1-z_1$ respectively. Moreover, they are disjoint. Since both $1+x_1-y_1$ and $1+x_1-z_1$ are not larger than $x_1$, by Lemma 2.3 it follows that
$$T_8\cup T_9\subset D'_1(C,X').$$ Then, applying (4.14) and the assumption on ${\bf x}_1$ we get
\begin{eqnarray*}
\hspace{1.5cm}{\rm vol}(D'_1(C,X'))&\hspace{-0.2cm}\ge &\hspace{-0.2cm} {\rm vol}(T_2\cup T_3\cup T_8\cup T_9)\\
&\hspace{-0.2cm}=&\hspace{-0.2cm}{1\over 6}\left((y_1-1)^3+(z_1-1)^3-(1-x_1)^3+(1+x_1-y_1)^3+(1+x_1-z_1)^3\right)\\
&\hspace{-0.2cm}\ge &\hspace{-0.2cm} {1\over 6}\left(2\left(1-{x_1\over 2}\right)^3-(1-x_1)^3+2\left({{3x_1}\over 2}-1\right)^3\right)\\
&\hspace{-0.2cm}\ge &\hspace{-0.2cm} {1\over 6} \left({5\over 9}-{4\over 9}\sqrt{1\over {10}}\right),\hspace{8.4cm} (4.16)
\end{eqnarray*}
where all the equalities hold if and only if
${\bf x}_1=\left({2\over 3}+{2\over 3}\sqrt{1\over {10}},\ {5\over 3}-{1\over 3}\sqrt{1\over {10}},\ {5\over 3}-{1\over 3}\sqrt{1\over {10}}\right)$ and ${\bf x}_2=\left(2,\ {1\over 3}-{2\over 3}\sqrt{1\over {10}},\right.$ $\left.{1\over 3}-{2\over 3}\sqrt{1\over {10}}\right).$

\medskip\noindent
{\bf Subcase 2.2.} $y_2\le z_1-2.$ Then the assumption ${\rm int}(F_3+{\bf x}_1)\cap J\not=\emptyset $ implies
$$y_1-2\le y_2\le z_1-2.\eqno (4.17)$$
Let $T_0'$ to be the orthogonal tetrahedron with vertices $(1,y_2+1,z_1-1)$, $(2, y_2+1, z_1-1)$, $(1,y_2+2,z_1-1)$ and $(1, y_2+1, z_1-2)$, let $S_2$, $S_8$, $S_9$, $S_{10}$ and $S_{11}$ denote the halfspaces $\{ (x,y,z):\ z\ge 0\}$, $\{ (x,y,z):\ x\le x_1+1\}$, $\{ (x,y,z):\ z\ge z_2-1\}$, $\{ (x,y,z):\ x-y-z\ge 1\}$ and $\{ (x,y,z):\ x-y+z\ge 2\}$, respectively, and define
$$P_3=P_1\cap T_0'\cap S_8\cap S_9.$$
It can be verified by Lemma 2.3 and Lemma 2.4 that
$$P_3\subseteq D'_1(C,X').$$

By $x_2-y_2+z_2\le 4$, $(x_2-x_1)-(y_2-y_1)-(z_2-z_1)=4$ and $x_1+y_1+z_1=4$ one can deduce $z_2\le 2-x_1$.
Then by routine computations based on (4.17) it can be shown that
\begin{eqnarray*}
{\rm vol}(T'_0\setminus \{S_2\cap S_9\})&\hspace{-0.2cm}\le &\hspace{-0.2cm}\max \left\{ {1\over 6}(2-z_1)^3,\  {1\over 6}(1-z_1+z_2)^3\right\}
\le {1\over 6}(y_1-1)^3,\\
{\rm vol}(T'_0\cap S_8\cap S_{10})&\hspace{-0.2cm}\le &\hspace{-0.2cm}{1\over 2}\left({{1-z_1-y_2}\over {\sqrt{2}}}\right)^2x_1\le {1\over 6}\times {3\over 2}(2z_1-3)^2x_1,\\
{\rm vol}(T'_0\cap S_{11})&\hspace{-0.2cm}\le &\hspace{-0.2cm}{1\over 6}\times 2\left({{z_1-y_2-2}\over 2}\right)^3\le {1\over 6}\times {1\over 4}(z_1-y_1)^3,
\end{eqnarray*}
and therefore
$${\rm vol}(P_3)\ge {1\over 6}\left(1-(1-x_1)^3-(1+y_2)^3-(y_1-1)^3-{3\over 2}(2z_1-3)^2x_1-{1\over 4}(z_1-y_1)^3\right).$$
Thus, by (4.14), (4.17) and the assumptions on ${\bf x}_1$ and ${\bf x}_2$ we obtain
\begin{eqnarray*}
\hspace{0.9cm}{\rm vol}(D'_1(C,X'))&\hspace{-0.2cm}\ge &\hspace{-0.2cm}{\rm vol}(T_2\cup T_3\cup P_3)\\
&\hspace{-0.2cm}\ge &\hspace{-0.2cm}{1\over 6}\left( 1+(z_1-1)^3-2(1-x_1)^3-(y_2+1)^3-{3\over 2}(2z_1-3)^2x_1-{1\over 4}(z_1-y_1)^3\right)\\
&\hspace{-0.2cm}\ge &\hspace{-0.2cm} {1\over 6} \left(1-2 (1-x_1)^3-{3\over 2}(2z_1-3)^2x_1-{1\over 4}(2\beta -3)^3\right)\\
&\hspace{-0.2cm}\ge &\hspace{-0.2cm}{1\over 6}\left(1-{3\over 2}(2\beta -3)^2-{1\over 4}(2\beta -3)^3\right)\\
&\hspace{-0.2cm}> &\hspace{-0.2cm}{1\over 6}\left({5\over 9}-{4\over 9}\sqrt{1\over {10}}\right).\hspace{9cm} (4.18)
\end{eqnarray*}

As a conclusion of (4.16) and (4.18), in this case we have
$${\rm vol}(D'_1(C,X'))\ge {1\over 6}\left({5\over 9}-{4\over 9}\sqrt{1\over {10}}\right).$$

\noindent
{\bf Case 3.} ${\rm int}(F_0)\cap J\not=\emptyset $. Assume that $C+{\bf x}_2$, where ${\bf x}_2=(2,y_2,z_2)\in X'$, touches $C$ at some interior points of $F_0$. We divide the region $\{ (y_2, z_2): \ -1\le y_2\le 2,\ -1\le z_2\le 2\}$ into nine parts as illustrated by Figure 6 and consider the following corresponding subcases.

\begin{figure}[hc]
\includegraphics[height=5.8cm,width=5.5cm,angle=0]{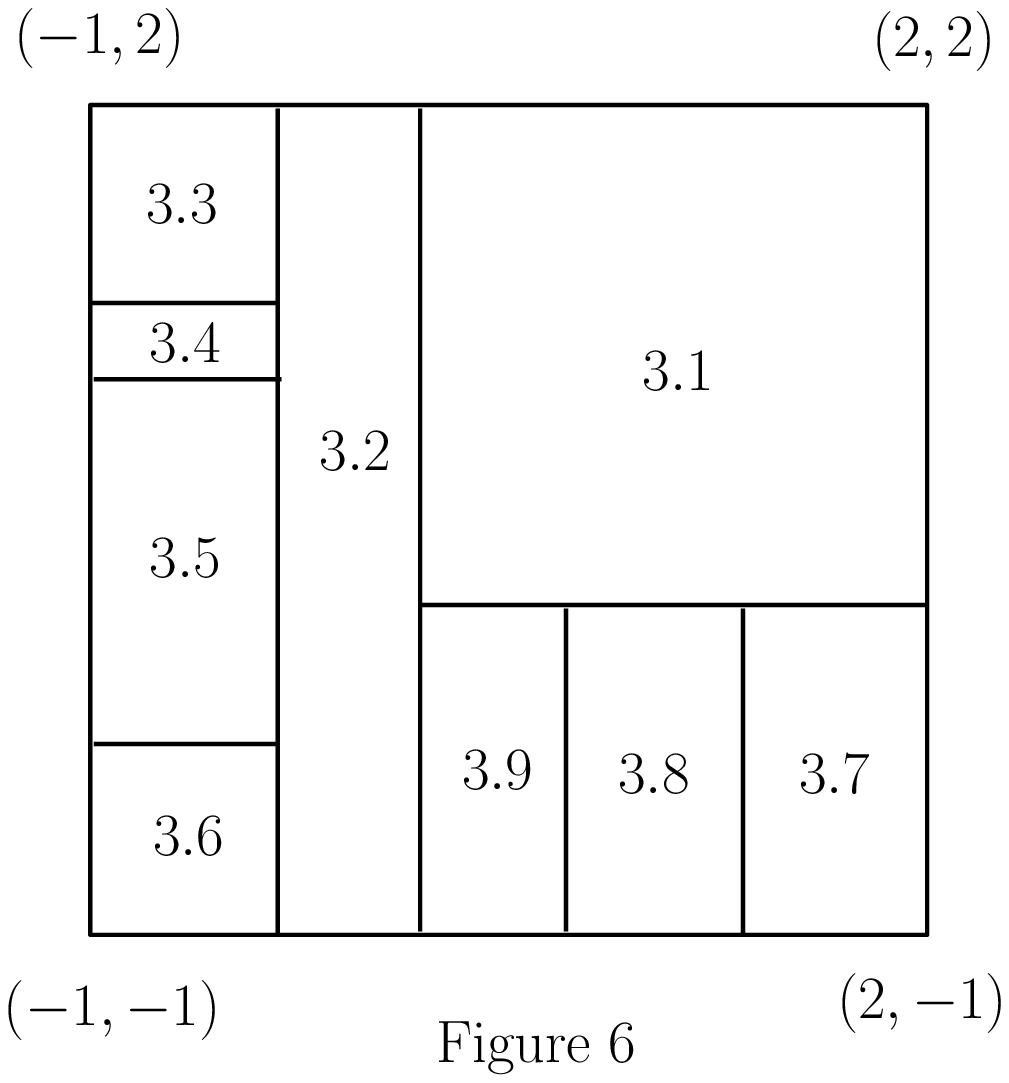}
\end{figure}

\medskip\noindent
{\bf Subcase 3.1.} {\it $y_2\ge z_1-2$ and $z_2\ge z_1-2$.} Then, the two translates $C+{\bf x}_1$ and $C+{\bf x}_2$ are separated by the hyperplane $(x-x_1)-(y-y_1)-(z-z_1)=2$. Let $T_8$ and $T_9$ be the tetrahedra defined in Subcase 2.1, we also have
$$T_8\cup T_9\subset D'_1(C,X').$$
If, on the contrary, there is a translate $C+{\bf x}$ satisfying both
$$(C+{\bf x})\cap {\rm int}(T_8)\not= \emptyset $$
and
$$(C+{\bf x})\cap ({\rm int}(C)+{\bf x}_i)=\emptyset, \qquad i=0,\ 1,\ 2,$$
by reducing the $x$-coordinate of ${\bf x}$ one can assume that $C+{\bf x}$ touches $(C+{\bf x}_0)\cup (C+{\bf x}_1)\cup (C+{\bf x}_2)$ at its boundary. By considering cases with respect to the facet of the touching points, one can reach contradictions one by one. For example, if $C+{\bf x}$ touches $C+{\bf x}_2$ at some interior point of $F_1+{\bf x}_2$, $x<x_1+2$, and ${\bf x}_1{\bf x}\prec F_5+{\bf x}_1$, we get
$$\left\{
\begin{array}{ll}
(x-2)+(y-y_2)+(z-z_2)=4, &\\
x<x_1+2,&\\
z\le z_1-2,
\end{array}
\right.$$
and therefore
$$z_2=x+y+z-y_2-6<(x_1+z_1)+(y-y_2)-6<-1,$$
which contradicts the assumption that ${\bf x}_2\in P_2$. In fact, this example is the only nontrivial case.

Similar to Subcase 2.1 we obtain
\begin{eqnarray*}
\hspace{1.3cm}{\rm vol}(D'_1(C,X'))&\hspace{-0.2cm}\ge &\hspace{-0.2cm} {\rm vol}(T_2\cup T_3\cup T_8\cup T_9)\\
&\hspace{-0.2cm}=&\hspace{-0.2cm}{1\over 6}\left((y_1-1)^3+(z_1-1)^3-(1-x_1)^3+(1+x_1-y_1)^3+(1+x_1-z_1)^3\right)\\
&\hspace{-0.2cm}\ge &\hspace{-0.2cm} {1\over 6} \left({5\over 9}-{4\over 9}\sqrt{1\over {10}}\right),\hspace{8.6cm} (4.19)\end{eqnarray*}
where all the equalities hold if and only if ${\bf x}_1=\left({2\over 3}+{2\over 3}\sqrt{1\over {10}},\ {5\over 3}-{1\over 3}\sqrt{1\over {10}},\ {5\over 3}-{1\over 3}\sqrt{1\over {10}}\right)$ and ${\bf x}_2=\left(2,\ {1\over 3}-{2\over 3}\sqrt{1\over {10}},\right.$ $\left.{1\over 3}-{2\over 3}\sqrt{1\over {10}}\right).$

\medskip\noindent
{\bf Subcase 3.2.} $y_1-2\le y_2\le z_1-2.$ Let $P_3$, $T_0'$, $S_2$, $S_8$, $S_9$, $S_{10}$ and $S_{11}$ be the polytopes defined in Subcase 2.2, here we also have
$$P_3\subseteq D'_1(C,X').$$

If, on the contrary, there is a translate $C+{\bf x}$ satisfying both
$$(C+{\bf x})\cap {\rm int}(P_3)\not= \emptyset $$
and
$$(C+{\bf x})\cap ({\rm int}(C)+{\bf x}_i)=\emptyset, \qquad i=0,\ 1,\ 2,$$
by reducing the $x$-coordinate of ${\bf x}$ one can assume that $C+{\bf x}$ touches $(C+{\bf x}_0)\cup (C+{\bf x}_1)\cup (C+{\bf x}_2)$ at its boundary. By considering cases with respect to the facet of the touching points, one can reach contradictions one by one. For example, if $C+{\bf x}$ touches $C+{\bf x}_2$ at some interior point of $F_1+{\bf x}_2$, $x<x_1+2$, and ${\bf x}_1{\bf x}\prec F_3+{\bf x}_1$, we have
$$\left\{
\begin{array}{ll}
(x-2)+(y-y_2)+(z-z_2)=4, &\\
x<x_1+2,&\\
(x-x_1)-(y-y_1)-(z-z_1)\ge 4,
\end{array}
\right.$$
and therefore
$$z_2=x+(y+z)-y_2-6\le 2x-x_1+y_1+z_1-y_2-10<4+x_1+y_1+z_1-y_1+2-10<-y_1\le -1,$$
which contradicts the assumption that ${\bf x}_2\in P_2$.

If $z_1-z_2>2$ and $2-x_1-(y_2+1-y_1)-(z_2+1-z_1)<2,$ it can be deduced that $y_1-x_1<y_2\le z_1-2$ and $y_1<1$, which contradicts the assumption ${\bf x}_1\in \triangle_4\cap H$. Thus we have ($z_1-z_2\le 2$ and $y_1-y_2\le 2$ together implies the same conclusion)
$$2-x_1-(y_2+1-y_1)-(z_2+1-z_1)\ge 2,$$
$$1-z_1+z_2\le y_1-x_1-y_2-1\le 1-x_1$$
and therefore
$${\rm vol}(T'_0\setminus \{S_2\cap S_9\})\le \max \left\{ {1\over 6}(2-z_1)^3,\  {1\over 6}(1-z_1+z_2)^3\right\}
\le {1\over 6}\left((2-z_1)^3+(1-x_1)^3\right).$$
In addition, similar to Subcase 2.2, we have
$${\rm vol}(T'_0\cap S_8\cap S_{10})\le {1\over 2}\left({{1-z_1-y_2}\over {\sqrt{2}}}\right)^2x_1\le {1\over 6}\times {3\over 2}(2z_1-3)^2x_1,$$
$${\rm vol}(T'_0\cap S_{11})\le {1\over 6}\times 2\left({{z_1-y_2-2}\over 2}\right)^3\le {1\over 6}\times {1\over 4}(z_1-y_1)^3$$
and therefore
\begin{eqnarray*}
{\rm vol}(P_3)&\hspace{-0.2cm}\ge &\hspace{-0.2cm}{1\over 6}\left(1-2(1-x_1)^3-(1+y_2)^3-(2-z_1)^3-{3\over 2}(2z_1-3)^2x_1-{1\over 4}(z_1-y_1)^3\right)\\
&\hspace{-0.2cm}\ge &\hspace{-0.2cm}{1\over 6}\left(1-2(1-x_1)^3-(z_1-1)^3-(2-z_1)^3-{3\over 2}(2z_1-3)^2x_1-{1\over 4}(z_1-y_1)^3\right).
\end{eqnarray*}

Thus, by (4.14) and the assumptions on ${\bf x}_1$ we obtain
\begin{eqnarray*}
\hspace{0.6cm}{\rm vol}(D'_1(C,X'))&\hspace{-0.2cm}\ge &\hspace{-0.2cm}{\rm vol}(T_2\cup T_3\cup P_3)\\
&\hspace{-0.2cm}\ge &\hspace{-0.2cm}{1\over 6}\left( 1+(y_1-1)^3-3(1-x_1)^3-(2-z_1)^3-{3\over 2}(2z_1-3)^2x_1-{1\over 4}(z_1-y_1)^3\right)\\
&\hspace{-0.2cm}\ge &\hspace{-0.2cm} {1\over 6} \left(1-3(1-x_1)^3-{3\over 2}(2z_1-3)^2x_1-{1\over 4}(z_1-y_1)^3\right)\\
&\hspace{-0.2cm}\ge &\hspace{-0.2cm}{1\over 6}\left(1-{3\over 2}(2\beta -3)^2-{1\over 4}(2\beta -3)^3\right)\\
&\hspace{-0.2cm}> &\hspace{-0.2cm}{1\over 6}\left({5\over 9}-{4\over 9}\sqrt{1\over {10}}\right).\hspace{9cm} (4.20)
\end{eqnarray*}

\medskip\noindent
{\bf Subcase 3.3.} {\it $-1\le y_2\le y_1-2$ and $z_2\ge 1.25$}. By the assumption ${\bf x}_2=(2,y_2,z_2)\in P_2$ we get
$$2-y_2+z_2< 4.$$
Thus one can deduce
$$1.25\le z_2< y_2+2\le y_1,$$
which implies
$$-0.75\le y_2$$
and
$$-0.5\le y_2+z_2-1\le 2y_1-3\le 1-\alpha .$$
Let $T_0^*$ denote the orthogonal tetrahedron with vertices $(1,y_2+1, z_2-1)$, $(2, y_2+1, z_2-1)$, $(1,y_2,z_2-1)$ and $(1,y_2+1,z_2)$, let $P_1$ be the polytope defined above (4.11), let $S_8$ be the halfspace defined in Subcase 2.2, and define
$$P_4=P_1\cap T_0^*\cap S_8.$$
By Lemma 2.3 and Lemma 2.4 it can be verified that
$$P_4\subset D'_1(C,X').$$

We recall $S_8=\{ (x,y,z):\ x\le x_1+1\}$, $S_{10}=\{ (x,y,z):\ x-y-z\ge 1\}$, and define $S'_6=\{ (x,y,z):\ x+y-z\ge 2\}$. By routine computations we have
$${\rm vol}(T^*_0\cap S'_6)={1\over 6}\times {1\over 4}(2+y_2-z_2)^3,$$
$${\rm vol}(T^*_0\cap S_8\cap S_{10})\le {1\over 6}\times {3\over 2}(y_2+z_2-1)^2x_1,$$
and
$${\rm vol}(P_4) \ge {1\over 6}\left(1-(z_2-1)^3-{1\over 4}(2+y_2-z_2)^3+y_2^3-{3\over 2}(y_2+z_2-1)^2x_1-(1-x_1)^3\right).$$
Then, by the assumptions on ${\bf x}_1$ and ${\bf x}_2$ and by (4.14), we get
\begin{eqnarray*}
\hspace{0.4cm}{\rm vol}(D'_1(C,X'))&\hspace{-0.2cm}\ge &\hspace{-0.2cm}{\rm vol}(T_2\cup T_3\cup P_4)\\
&\hspace{-0.2cm}\ge &\hspace{-0.2cm}{1\over 6}\left(1+(y_1-1)^3+(z_1-1)^3-2(1-x_1)^3-(z_2-1)^3-{1\over 4}(2+y_2-z_2)^3+y_2^3\right.\\
&\hspace{-0.2cm}&\hspace{-0.2cm}\left. -{3\over 2}(y_2+z_2-1)^2x_1\right)\\
&\hspace{-0.2cm}\ge &\hspace{-0.2cm}{1\over 6} \left(1+(y_1-1)^3+(z_1-1)^3-2(1-x_1)^3-(y_2+1)^3+y_2^3-{3\over 2}(y_2+z_2-1)^2\right)\\
&\hspace{-0.2cm}\ge &\hspace{-0.2cm}{1\over 6} \left(1+ (y_1-1)^3+(z_1-1)^3-2(1-x_1)^3-\left({1\over 4}\right)^3-
\left({3\over 4}\right)^3-{3\over 2}(y_2+z_2-1)^2\right)\\
&\hspace{-0.2cm}\ge &\hspace{-0.2cm} {1\over 6}\left(1+{1\over 4}-(1-\alpha)^3-\left({1\over 4}\right)^3-\left({3\over 4}\right)^3
-{3\over 2}\left({1\over 2}\right)^2\right)\\
&\hspace{-0.2cm}> &\hspace{-0.2cm}{1\over 6}\left({5\over 9}-{4\over 9}\sqrt{1\over {10}}\right). \hspace{9.6cm}(4.21)
\end{eqnarray*}

\medskip
\noindent
{\bf Subcase 3.4.} {\it $-1\le y_2\le y_1-2$ and $1\le z_2\le 1.25$}. Then by $(2, y_2, z_2)\in P_2$ one can deduce
$$1\le z_2\le \min \{ 1.25,\ y_2+2\}.$$

Let $P_3$ be the polytope defined in Subcase 2.2 and let $P_4$ be the polytope defined in Subcase 3.3. It can be shown that
$${\rm int}(P_3)\cap {\rm int}(P_4)=\emptyset $$
and, by Lemma 2.3 and Lemma 2.4,
$$P_i\subset D'_1(C,X'),\qquad i=3,\ 4.$$

By routine computations it can be shown that
$${\rm vol}(P_3)\ge {1\over 6}\left(1-(y_2+1)^3-{1\over 4}(z_1-y_2-2)^3-(1-z_1+z_2)^3-{3\over 2}(1-y_2-z_1)^2(z_1-z_2)-(1-x_1)^3\right)$$
and
$${\rm vol}(P_4)\ge {1\over 6}\left(1-(z_2-1)^3-{1\over 4}(2+y_2-z_2)^3+y_2^3-{3\over 2}(y_2+z_2-1)^2(y_2+1)-(1-x_1)^3\right).$$

If $y_2\ge -0.8$, we have
$$-0.8\le y_2\le y_1-2\le -{\alpha \over 2},$$
$$2+y_2-z_2\le 1+y_2\le 1-{\alpha \over 2}$$
and
$$-0.8\le y_2+z_2-1\le {1\over 4}-{\alpha \over 2}.$$
Then, together with (4.14), we get
\begin{eqnarray*}
\hspace{1.3cm}{\rm vol}(D'_1(C,X'))&\hspace{-0.2cm}\ge &\hspace{-0.2cm}{\rm vol}(T_2\cup T_3\cup P_4)\\
&\hspace{-0.2cm}\ge &\hspace{-0.2cm} {1\over 6}\left({5\over 4}-\left({1\over 4}\right)^3-{1\over 4}\left(1-{\alpha \over 2}\right)^3
+y_2^3-{3\over 2}y_2^2(y_2+1)-(1-\alpha )^3\right)\\
&\hspace{-0.2cm}\ge &\hspace{-0.2cm} {1\over 6}\left({5\over 4}-\left({1\over 4}\right)^3-{1\over 4}\left(1-{\alpha \over 2}\right)^3
-\left({4\over 5}\right)^3-{3\over 2}\left({4\over 5}\right)^2{1\over 5}-(1-\alpha )^3\right)\\
&\hspace{-0.2cm}>&\hspace{-0.2cm} {1\over 6}\left({5\over 9}-{4\over 9}\sqrt{1\over {10}}\right). \hspace{8.7cm}(4.22)
\end{eqnarray*}

When $-1\le y_2\le -0.8$ and $1\le z_2\le 1.25$, we have
$$z_1-y_2-2\le \beta -1$$
and
$$0\le 1-y_2-z_1\le 2-z_1\le {1\over 2}.$$
Then, together with (4.14), we get
\begin{eqnarray*}
\hspace{0.4cm}{\rm vol}(D'_1(C,X'))&\hspace{-0.2cm}\ge &\hspace{-0.2cm}{\rm vol}(T_2\cup T_3\cup P_3)\\
&\hspace{-0.2cm}\ge &\hspace{-0.2cm} {1\over 6}\left({5\over 4}-\left({1\over 5}\right)^3-{1\over 4}\left(\beta -1\right)^3
-(1-(z_1-z_2))^3-{3\over 2}\left({1\over 2}\right)^2(z_1-z_2)-(1-\alpha )^3\right)\\
&\hspace{-0.2cm}\ge &\hspace{-0.2cm} {1\over 6}\left({5\over 4}-\left({1\over 5}\right)^3-{1\over 4}\left(\beta -1\right)^3
-\left({3\over 4}\right)^3-{3\over 2}\left({1\over 2}\right)^2{1\over 4}-(1-\alpha )^3\right)\\
&\hspace{-0.2cm}>&\hspace{-0.2cm} {1\over 6}\times {1\over 2}. \hspace{11cm}(4.23)
\end{eqnarray*}

As a conclusion of (4.22) and (4.23), in this subcase we have
$${\rm vol}(D'_1(C, X'))> {1\over 6}\left({5\over 9}-{4\over 9}\sqrt{1\over {10}}\right).$$

\medskip\noindent
{\bf Subcase 3.5.} {\it $-1\le y_2\le y_1-2$ and $y_1-2\le z_2\le 1$}. Then we have
$$P_3\subset D'_1(C,X'),$$
where $P_3$ was defined in Subcase 2.2. However, this time we have
$$z_1-z_2\ge z_1-1$$
and
$${\rm vol}(P_3)\ge {1\over 6}\left(1-(1+y_2)^3-{1\over 4}(z_1-y_2-2)^3-(2-z_1)^3-
{3\over 4}(1-y_2-z_1)^2(4+y_2-z_1)-(1-x_1)^3\right).$$

By routine computations one can deduce
$$(1+y_2)^3+{1\over 4} (z_1-y_2-2)^3\le \max \left\{ {1\over 4} (z_1-1)^3,\ (z_1-1)^3\right\}\le (z_1-1)^3$$
and
$$(1-y_2-z_1)^2(4+y_2-z_1)\le \max\{ (2-z_1)^2(3-z_1),\ (1-x_1)^2(2+y_1-z_1)\} \le \left({1\over 2}\right)^2
{3\over 2}.$$
Thus, we get
\begin{eqnarray*}
\hspace{1.4cm}{\rm vol}(D'_1(C,X'))&\hspace{-0.2cm}\ge &\hspace{-0.2cm} {\rm vol}(T_2\cup T_3\cup P_3)\\
&\hspace{-0.2cm}\ge &\hspace{-0.2cm}{1\over 6}\left(1+(y_1-1)^3+(z_1-1)^3-2(1-x_1)^3 -(1+y_2)^3-{1\over 4}(z_1-y_2-2)^3\right.\\
&\hspace{-0.2cm}&\hspace{-0.2cm}-\left. (2-z_1)^3-{3\over 4}(1-y_2-z_1)^2(4+y_2-z_1)\right)\\
&\hspace{-0.2cm}\ge &\hspace{-0.2cm} {1\over 6}\left( 1+(y_1-1)^3-2(1-x_1)^3-(2-z_1)^3-{3\over 4}\left( {1\over 2}\right)^2{3\over 2}\right)\\
&\hspace{-0.2cm}\ge &\hspace{-0.2cm} {1\over 6}\left( 1+(2-\beta )^3-2(1-\alpha )^3-\left({1\over 2}\right)^3-
{3\over 4}\left({1\over 2}\right)^2{3\over 2}\right)\\
&\hspace{-0.2cm}> &\hspace{-0.2cm} {1\over 6} \times {1\over 2}. \hspace{10.3cm}(4.24)
\end{eqnarray*}

\medskip\noindent
{\bf Subcase 3.6.} {\it $-1\le y_2\le y_1-2$ and $z_2\le y_1-2$}. First, by ${\bf x}_2\in P_2$ and ${\bf x}_1\in \triangle_4\cap H$ one can deduce
$$y_2+z_2\ge -1, \eqno (4.25)$$
$${\alpha \over 2}-1\le 1-y_1\le -1-z_2\le y_2\le y_1-2\le -{\alpha \over 2}\eqno (4.26)$$
and
$${\alpha \over 2}-1\le 1-y_1\le -1-y_2\le z_2\le y_1-2\le -{\alpha \over 2}.\eqno (4.27)$$

Let $T_0^\circ$ denote the tetrahedron with vertices $(1,y_2+1,z_2+1)$, $(2,y_2+1, z_2+1)$, $(1,y_2, z_2+1)$
and $(1, y_2+1, z_2)$, let $S_8$ and $S_{11}$ be the halfspaces defined in Subcase 2.2, and define
$$P_5=P_1\cap T_0^\circ\cap S_8.$$
It can be verified by Lemma 2.3 and Lemma 2.4 that
$$P_5\subset D'_1(C,X').$$
By routine computations it can be deduced that
$${\rm vol}(T_0^\circ \cap S_8\cap S_{11})\le {1\over 2}\left({{y_2-z_2}\over {\sqrt{2}}}\right)^2x_1={1\over 4}(y_2-z_2)^2x_1$$
and therefore
$${\rm vol}(P_5)\ge {1\over 6} \left( 1+y_2^3+z^3_2-{3\over 2}(y_2-z_2)^2x_1-(1-x_1)^3\right).$$

Thus, by (4.14), (4.25), (4.26) and (4.27) we get
\begin{eqnarray*}
\hspace{1.2cm}{\rm vol}(D'_1(C, X'))&\hspace{-0.2cm}\ge &\hspace{-0.2cm} {\rm vol}(T_2\cup T_3\cup P_5)\\
&\hspace{-0.2cm}\ge &\hspace{-0.2cm} {1\over 6} \left( 1+(y_1-1)^3+(z_1-1)^3-2(1-x_1)^3+y_2^3+z^3_2-{3\over 2}(y_2-z_2)^2x_1\right) \\
&\hspace{-0.2cm}\ge &\hspace{-0.2cm} {1\over 6} \left( {5\over 4}-(1-\alpha)^3+y^3_2+z_2^3-{3\over 2}(1-\alpha )^2\right) \\
&\hspace{-0.2cm}\ge &\hspace{-0.2cm} {1\over 6} \left( {5\over 4}-(1-\alpha)^3 -\left({\alpha\over 2}\right)^3-\left(1-{\alpha\over 2}\right)^3-{3\over 2}(1-\alpha )^2\right) \\
&\hspace{-0.2cm}>&\hspace{-0.2cm} {1\over 6}\times {4\over 5}.\hspace{10.4cm}(4.28)
\end{eqnarray*}

\medskip\noindent
{\bf Subcase 3.7.} {\it $-1\le z_2\le z_1-2$ and $y_2\ge 1.2$}. By the assumption ${\bf x}_2=(2,y_2,z_2)\in P_2$ we get
$$2+y_2-z_2< 4.$$
Thus one can deduce
$$1.2\le y_2< z_2+2\le z_1,$$
which implies
$$-0.8\le z_2\le z_1-2$$
and
$$z_2+0.2\le y_2+z_2-1\le 2z_2+1\le 2z_1-3\le 2\beta -3.$$
Let $T_0^\bullet $ denote the orthogonal tetrahedron with vertices $(1,y_2-1, z_2+1)$, $(2, y_2-1, z_2+1)$, $(1,y_2,z_2+1)$ and $(1,y_2-1,z_2)$, let $P_1$ be the polytope defined above (4.11), let $S_8$ be the halfspace defined in Subcase 2.2, and define
$$P_6=P_1\cap T_0^\bullet \cap S_8.$$
By Lemma 2.3 and Lemma 2.4 it can be verified that
$$P_6\subset D'_1(C,X').$$

We recall $S_2=\{ (x,y,z):\ z\ge 0\}$, $S_{10}=\{ (x,y,z):\  x-y-z\ge 1\}$ and $S_{11}=\{ (x,y,z):\ x-y+z\ge 2\}$. By routine computations we have
$${\rm vol}(T^\bullet_0\cap S_{11})={1\over 6}\times {1\over 4}(2-y_2+z_2)^3,$$
$${\rm vol}(T^\bullet_0\cap S_2\cap S_{10})\le {1\over 6}\times {3\over 2}(y_2+z_2-1)^2(z_2+1),$$
$${\rm vol}(P_6) \ge {1\over 6}\left(1-(y_2-1)^3-{1\over 4}(2-y_2+z_2)^3+z_2^3-{3\over 2}(y_2+z_2-1)^2(z_2+1)-(1-x_1)^3\right)$$
and therefore
\begin{eqnarray*}
{\rm vol}(D'_1(C,X'))&\hspace{-0.2cm}\ge &\hspace{-0.2cm}{\rm vol}(T_2\cup T_3\cup P_6)\\
&\hspace{-0.2cm}\ge &\hspace{-0.2cm}{1\over 6}\left(1+(y_1-1)^3+(z_1-1)^3-2(1-x_1)^3-(y_2-1)^3-{1\over 4}(2-y_2+z_2)^3+z_2^3\right.\\
&\hspace{-0.2cm}&\hspace{-0.2cm}\left. -{3\over 2}(y_2+z_2-1)^2(z_2+1)\right)\\
&\hspace{-0.2cm}\ge &\hspace{-0.2cm} {1\over 6} \left(1+(y_1-1)^3+(z_1-1)^3-2(1-x_1)^3-(z_2+1)^3+z_2^3-{3\over 2}(y_2+z_2-1)^2(z_2+1)\right).
\end{eqnarray*}

When $-0.8\le z_2\le -0.4$, we have
$$z_2+0.2\le y_2+z_2-1\le 2z_2+1\le -(z_2+0.2)$$
and
\begin{eqnarray*}
\hspace{0.3cm}{\rm vol}(D'_1(C, X'))&\hspace{-0.2cm}\ge &\hspace{-0.2cm} {1\over 6} \left(1+(y_1-1)^3+(z_1-1)^3-2(1-x_1)^3-(z_2+1)^3+z_2^3-{3\over 2}(z_2+0.2)^2(z_2+1)\right)\\
&\hspace{-0.2cm}\ge &\hspace{-0.2cm} {1\over 6}\left( {5\over 4} -(1-\alpha)^3-\left({1\over 5}\right)^3-\left({4\over 5}\right)^3
-{3\over 2}\left({3\over 5}\right)^2{1\over 5}\right)\\
&\hspace{-0.2cm}\ge &\hspace{-0.2cm}{1\over 6}\times {3\over 5}. \hspace{11.1cm}(4.29)
\end{eqnarray*}
When $-0.4\le z_2\le z_1-2$, we get
$$-(2z_2+1)\le z_2+0.2\le y_2+z_2-1\le 2z_2+1$$
and
\begin{eqnarray*}
\hspace{0.3cm}{\rm vol}(D'_1(C,X'))&\hspace{-0.2cm}\ge &\hspace{-0.2cm} {1\over 6} \left(1+ (y_1-1)^3+(z_1-1)^3-2(1-x_1)^3-(z_2+1)^3+z_2^3-{3\over 2}(2z_2+1)^2(z_2+1)\right)\\
&\hspace{-0.2cm}\ge &\hspace{-0.2cm} {1\over 6} \left(1+ (y_1-1)^3-2(1-x_1)^3-(2-z_1)^3-{3\over 2}(2z_1-3)^2(z_1-1)\right)\\
&\hspace{-0.2cm}\ge &\hspace{-0.2cm} {1\over 6}\left(1+(2-\beta )^3-2(1-\alpha)^3-\left(2-\beta \right)^3-{3\over 2}\left(2\beta -3\right)^2(\beta -1)\right)\\
&\hspace{-0.2cm}> &\hspace{-0.2cm}{1\over 6}\times {1\over 2}. \hspace{11.1cm}(4.30)
\end{eqnarray*}

As a conclusion of (4.29) and (4.30), in this subcase we have
$${\rm vol}(D'_1(C, X'))> {1\over 6}\times {1\over 2}.$$

\medskip\noindent
{\bf Subcase 3.8.} $-1\le z_2\le z_1-2$ {\it and} $0\le y_2\le 1.2$. We define
\begin{eqnarray*}
P_7&\hspace{-0.2cm}=&\hspace{-0.2cm}P_1\cap \left\{ (x,y,z):\ 1\le x\le x_1+1,\ \max\{ 0, y_2-1\} \le y\le 1, \ z_2+1\le z\le z_1-1,\right.\\
&\hspace{-0.2cm}&\hspace{-0.2cm}\left. (x-1)-(y-y_1+1)+(z-z_2-1)\le 1\right\},\\
P_8&\hspace{-0.2cm}=&\hspace{-0.2cm}P_1\cap \left\{ (x,y,z):\ 1\le x\le x_1+1,\ \max\{ 0, y_2-1\}\le y\le 1,\ 0\le z\le z_2+1\right\}\setminus \{ C+{\bf x}_2\},\\
\overline{P_7}&\hspace{-0.2cm}=&\hspace{-0.2cm}P_7\cap \{ (x,y,z):\ z=z_2+1\}\qquad {\rm and}\qquad \overline{P_8}=P_8\cap \{ (x,y,z):\ z=z_2+1\},
\end{eqnarray*}
and let $s(\overline{P_i})$ denote the area of $\overline{P_i}$. Clearly, we have
$${\rm int}(P_7)\cap {\rm int}(P_8)= \emptyset .$$
By Lemmas 2.3 and 2.4 it can be shown that
$$P_i\subseteq D'_1(C, X')$$
holds for both $i=7$ and $8$.

When $z_2\ge y_1-2$, we have
$$P_7=P_1\cap \left\{ (x,y,z):\ 1\le x\le x_1+1,\ \max\{ 0,\ y_2-1\} \le y\le 1, \ z_2+1\le z\le z_1-1\right\}$$
and
$$\overline{P_8}\subseteq \overline{P_7}.$$
Now we consider ${\rm vol}(P_7)$ and ${\rm vol}(P_8)$ as functions of $x_1$, $y_1$, $z_1$, $y_2$ and $z_2$. It can be shown that
$${{\partial {\rm vol}(P_7)}\over {\partial z_2}}=-s(\overline{P_7}),$$
$${{\partial {\rm vol}(P_8)}\over {\partial z_2}}\le s(\overline{P_8}),$$
$${{\partial ({\rm vol}(P_7)+{\rm vol}(P_8))}\over {\partial z_2}}\le -s(\overline{P_7})+s(\overline{P_8})\le 0$$
and therefore
$${\rm vol}(P_7)+{\rm vol}(P_8)\ge \left({\rm vol}(P_7)+{\rm vol}(P_8)\right)\bigg|_{z_2=z_1-2} \ge {\rm vol}(P_8)\bigg|_{z_2=z_1-2}.$$
When $0\le y_2\le 1$, based on a figure similar to Figure 5 one can deduce
$${\rm vol}(P_8)\bigg|_{z_2=z_1-2}\ge {1\over 6}\times 2\left(1-{{z_1}\over 3}\right)^3.$$
When $1\le y_2\le 1.3$, we get
$${\rm vol}(P_8)\bigg|_{z_2=z_1-2}\ge {1\over 6}\times (0.7^3-3\times 0.2^3)\ge {1\over 6}\times 2\left(1-{{z_1}\over 3}\right)^3,$$
since $z_1\ge 1.5$. Thus, we have
$${\rm vol}(P_7)+{\rm vol}(P_8)\ge {1\over 6}\times 2\left(1-{{z_1}\over 3}\right)^3.$$
Then, together with (4.14) we get
\begin{eqnarray*}
\hspace{1.7cm}{\rm vol}(D'_1(C, X'))&\hspace{-0.2cm}\ge &\hspace{-0.2cm}{\rm vol}(T_2\cup T_3\cup P_7\cup P_8)\\
&\hspace{-0.2cm}\ge &\hspace{-0.2cm}{1\over 6}\left((y_1-1)^3+(z_1-1)^3-(1-x_1)^3+2\left(1-{{z_1}\over 3}\right)^3\right)\\
&\hspace{-0.2cm}>&\hspace{-0.2cm}{1\over 6}\left({5\over 9}-{4\over 9}\sqrt{1\over {10}}\right). \hspace{8.3cm}(4.31)
\end{eqnarray*}

When $z_2\le y_1-2$ and $1\le y_2\le 1.2$. It can be shown that $P_7$ contains $(1,1, z_2+1)$, $(1, y_2-1, z_2+1)$,
$(1, y_2-1, z_1-1)$, $(1,1, z_1-1)$ and $(2.5+z_2, 0.5, z_2+1)$, and therefore
$${\rm vol}(P_7)\ge {1\over 3} (2-y_2)(z_1-z_2-2)(1.5+z_2).$$
If $-1\le z_2\le -0.8$, by (4.14) we get
\begin{eqnarray*}
\hspace{1.6cm}{\rm vol}(D'_1(C, X'))&\hspace{-0.2cm}\ge &\hspace{-0.2cm}{\rm vol}(T_2\cup T_3\cup P_7)\\
&\hspace{-0.2cm}\ge &\hspace{-0.2cm}{1\over 6}\left((y_1-1)^3+(z_1-1)^3-(1-x_1)^3+2(2-y_2)(z_1-z_2-2)(1.5+z_2)\right)\\
&\hspace{-0.2cm}\ge &\hspace{-0.2cm}{1\over 6}\left({1\over 4}+2(2-1.2)(1.5+0.8-2)(1.5-1)\right)\\
&\hspace{-0.2cm}>&\hspace{-0.2cm}{1\over 6}\left({5\over 9}-{4\over 9}\sqrt{1\over {10}}\right). \hspace{8.4cm}(4.32)
\end{eqnarray*}
If $-0.8\le z_2\le y_1-2\le -0.5$, let $P_6$ be the polytope defined in Subcase 3.7, by (4.14) we get
\begin{eqnarray*}
\hspace{1.6cm}{\rm vol}(D'_1(C, X'))&\hspace{-0.2cm}\ge &\hspace{-0.2cm}{\rm vol}(T_2\cup T_3\cup P_6)\\
&\hspace{-0.2cm}\ge &\hspace{-0.2cm}{1\over 6}\left({5\over 4} -(z_2+1)^3+z_2^3-{3\over 2}z_2^2(z_2+1)-(1-x_1)^3\right)\\
&\hspace{-0.2cm}\ge &\hspace{-0.2cm} {1\over 6}\left( {5\over 4} -\left({1\over 5}\right)^3-\left({4\over 5}\right)^3
-{3\over 2}\left({4\over 5}\right)^2{1\over 5}-(1-\alpha )^3\right)\\
&\hspace{-0.2cm}\ge &\hspace{-0.2cm}{1\over 6}\times {1\over 2}. \hspace{10.1cm}(4.33)
\end{eqnarray*}

When $z_2\le y_1-2$ and $0\le y_2\le 1$. Let $T_0^\star$ to be the tetrahedron with vertices $(1, y_1-1, z_2+1)$, $(2, y_1-1, z_2+1)$, $(1, y_1-2, z_2+1)$ and $(1, y_1-1, z_2+2)$, and define
$$P_9=P_1\cap T_0^\star \cap S_8.$$
By routine arguments it can be shown that
$$P_9\subseteq D'_1(C,X')$$
and
$${\rm vol}(P_9)\ge {1\over 6}\left(1-(2-y_1)^3-(z_2+1)^3-{1\over 4}(y_1-z_2-2)^3-{3\over 2}(1-y_1-z_2)^2(y_1-1)-(1-x_1)^3\right).$$
In addition, we have
$$(z_2+1)^3+{1\over 4}(y_1-z_2-2)^3\le (y_1-1)^3$$
and
$$y_1-2\le 3-2y_1\le 1-y_1-z_2\le 2-y_1.$$
Thus, combined with (4.14), we get
\begin{eqnarray*}
\hspace{1.5cm}{\rm vol}(D'_1(C, X'))&\hspace{-0.2cm}\ge &\hspace{-0.2cm}{\rm vol}(T_2\cup T_3\cup P_9)\\
&\hspace{-0.2cm}\ge &\hspace{-0.2cm}{1\over 6}\left(1+(z_1-1)^3-2(1-x_1)^3-(2-y_1)^3-{3\over 2}(2-y_1)^2(y_1-1)\right)\\
&\hspace{-0.2cm}>&\hspace{-0.2cm}{1\over 6}\times {7\over {10}}. \hspace{10cm}(4.34)
\end{eqnarray*}

As a conclusion of (4.31), (4.32), (4.33) and (4.34), in this subcase we have
$${\rm vol}(D'_1(C, X'))> {1\over 6}\left({5\over 9}-{4\over 9}\sqrt{1\over {10}}\right).$$

\medskip\noindent
{\bf Subcase 3.9.} $-1\le z_2\le z_1-2$ {\it and} $z_1-2\le y_2\le 0$. Let $P_5$ and $P_9$ be the polytopes defined in Subcases 3.6 and 3.8, respectively. In other words,
$$P_5=P_1\cap T^\circ_0\cap S_8$$
and
$$P_9=P_1\cap T_0^\star \cap S_8.$$

Clearly, we have
$$y_2+z_2>-1$$
and
$${\rm int}(P_5)\cap {\rm int}(P_9)= \emptyset .$$
By Lemmas 2.3 and 2.4 it can be shown that
$$P_i\subseteq D'_1(C, X')$$
holds for both $i=5$ and $9$.

When $-1\le z_2\le -0.8$, we have
$$z_2\le 1-\beta \le y_1-2, \qquad y_2\ge -1-z_2\ge -0.2,$$
$${\rm vol}(P_9)\ge {1\over 6}\left(1-(2-y_1)^3-(z_2+1)^3-{1\over 4}(y_1-z_2-2)^3-{3\over 2}(1-y_1-z_2)^2(y_1-1)-(1-x_1)^3\right)$$
and therefore, by (4.14),
\begin{eqnarray*}
\hspace{0.5cm}{\rm vol}(D'_1(C, X'))&\hspace{-0.2cm}\ge &\hspace{-0.2cm}{\rm vol}(T_2\cup T_3\cup P_9)\\
&\hspace{-0.2cm}\ge &\hspace{-0.2cm}{1\over 6}\left({5\over 4}-(2-y_1)^3-\left({1\over 5}\right)^3-{1\over 4}\left(1-{\alpha\over 2}\right)^3
-{3\over 2}(1-y_1-z_2)^2(y_1-1)-(1-\alpha )^3\right)\\
&\hspace{-0.2cm}\ge &\hspace{-0.2cm}{1\over 6}\left( {5\over 4}-(\beta -1)^3-\left({1\over 5}\right)^3-{1\over 4}\left(1-{\alpha\over 2}\right)^3-{3\over 2}(\beta -1)^2(2-\beta )-(1-\alpha )^3\right)\\
&\hspace{-0.2cm}> &\hspace{-0.2cm} {1\over 6} \left({5\over 9}-{4\over 9}\sqrt{1\over {10}}\right).\hspace{9.1cm}(4.35)
\end{eqnarray*}

If $-0.8\le z_2\le z_1-2$ and $z_1-2\le y_2\le 0$, we have
$$0\le y_2-z_2\le -z_2,$$
$${\rm vol}(P_5)\ge {1\over 6}\left(1+y_2^3+z_2^3-{3\over 2}(y_2-z_2)^2(z_2+1)-(1-x_1)^3\right)$$
and, together with (4.14),
\begin{eqnarray*}
\hspace{1.7cm}{\rm vol}(D'_1(C, X'))&\hspace{-0.2cm}\ge &\hspace{-0.2cm}{\rm vol}(T_2\cup T_3\cup P_5)\\
&\hspace{-0.2cm}\ge &\hspace{-0.2cm}{1\over 6}\left( {5\over 4}+y_2^3+z_2^3-{3\over 2}(y_2-z_2)^2(z_2+1)-(1-x_1)^3 \right)\\
&\hspace{-0.2cm}\ge &\hspace{-0.2cm}{1\over 6}\left( {5\over 4}-\left({1\over 5}\right)^3-\left({4\over 5}\right)^3
-{3\over 2}z_2^2(z_2+1)-(1-\alpha )^3\right)\\
&\hspace{-0.2cm}\ge &\hspace{-0.2cm}{1\over 6}\left( {5\over 4}-\left({1\over 5}\right)^3-\left({4\over 5}\right)^3
-{3\over 2}\left({2\over 3}\right)^2{1\over 3}-(1-\alpha )^3\right)\\
&\hspace{-0.2cm}> &\hspace{-0.2cm}{1\over 6}\left({5\over 9}-{4\over 9}\sqrt{1\over {10}}\right). \hspace{8cm}(4.36)
\end{eqnarray*}

By (4.35) and (4.36), in Subcase 3.9 we get
$${\rm vol}(D'_1(C, X'))> {1\over 6} \left({5\over 9}-{4\over 9}\sqrt{1\over {10}}\right).$$

\bigskip
As a conclusion of these subcases we get
$${\rm vol}(D'_1(C,X'))\ge {1\over 6} \left({5\over 9}-{4\over 9}\sqrt{1\over {10}}\right)$$
for Case 3, where the equality holds if and only if
${\bf x}_1=\left({2\over 3}+{2\over 3}\sqrt{1\over {10}},\ {5\over 3}-{1\over 3}\sqrt{1\over {10}},\
{5\over 3}-{1\over 3}\sqrt{1\over {10}}\right)$ and ${\bf x}_2=\left(2,\ {1\over 3}-{2\over 3}\sqrt{1\over {10}},\ {1\over 3}-{2\over 3}\sqrt{1\over {10}}\right).$
Lemma 4.4 is proved. \hfill{$\Box $}

\medskip
\noindent
{\bf Lemma 4.5.} {\it If ${\bf x}_1\in \triangle_3\cap H$, then we have}
$${\rm vol}(D_1(C,X))\ge {1\over 6} \left({5\over 9}-{4\over 9}\sqrt{1\over {10}}\right).$$

\medskip\noindent
{\bf Proof.} In this case $C+{\bf x}_1$ touches $C$ at an hexagon with vertices $(1, 2-z_1, z_1-1)$, $(1, y_1-1,
2-y_1)$, $(2-y_1, y_1-1, 1)$, $(x_1-1, 2-x_1, 1)$, $(x_1-1, 1,2-x_1)$ and $(2-z_1, 1, z_1-1)$, as shown in Figure 2. Then, by Lemma 2.3 and Lemma 2.4 it can be shown that $D'_1(C, X')$ contains four tetrahedra $T_1$, $T_2$, $T_3$ and $T_4$, where $T_1$ has vertices $(x_1-1,1,1)$, $(0,1,1)$, $(x_1-1,2-x_1,1)$ and $(x_1-1,1,2-x_1)$, $T_2$ has vertices $(1,y_1-1,1)$, $(2-y_1,y_1-1,1)$, $(1,0,1)$ and $(1,y_1-1,2-y_1)$, $T_3$ has vertices $(1,1,z_1-1)$, $(2-z_1,1,z_1-1)$, $(1,2-z_1,z_1-1)$ and $(1,1,0)$, and $T_4$ has vertices $(1,y_1-1,z_1-1)$, $(x_1,y_1-1,z_1-1)$, $(1,2-z_1,z_1-1)$ and $(1,y_1-1,2-y_1)$. Clearly, all $T_1$, $T_2$, $T_3$ and $T_4$ are homothetic to $T_0$ with ratios $x_1-1$, $y_1-1$, $z_1-1$ and $x_1-1$, respectively. Moreover, their interiors are pairwise disjoint. Thus, we have
$${\rm vol}(T_1\cup T_2\cup T_3\cup T_4)\ge {1\over 6}\left( 2(x_1-1)^3+(y_1-1)^3+(z_1-1)^3\right) \ge
{1\over 6}\times {{4\sqrt{2}+2}\over {25+22\sqrt{2}}},\eqno (4.37)$$
where the last equality holds if and only if ${\bf x}_1=\left({{2\sqrt{2}+2}\over {2\sqrt{2}+1}}, {{3\sqrt{2}+1}\over {2\sqrt{2}+1}}, {{3\sqrt{2}+1}\over {2\sqrt{2}+1}}\right).$ In particular, whenever $x_1\ge 1.585$ we have
$${\rm vol}(T_1\cup T_2\cup T_3\cup T_4)\ge {1\over 6}\left( 2(x_1-1)^3+(y_1-1)^3+(z_1-1)^3\right)>{1\over 6}
\left( {5\over 9}-{4\over 9}\sqrt{1\over {10}}\right).$$
Thus, in the rest of this proof we assume that
$$x_1\le 1.585.\eqno (4.38)$$

When ${\bf x}_1\in \triangle_3\cap H$, one can deduce
$$J\cap (F_0+{\bf x}_1)=\emptyset $$
and therefore by (4.13)
$$J\cap ({\rm int}(F_0)\cup ({\rm int}(F_3)+{\bf x}_1))\not= \emptyset .\eqno (4.39)$$
Now, we estimate ${\rm vol}(D'_1(C,X'))$ by considering two cases based on (4.39).

\medskip\noindent
{\bf Case 1.} $({\rm int}(F_3)+{\bf x}_1)\cap J\not= \emptyset .$ Assume that $C+{\bf x}_2$, where ${\bf x}_2=
(x_2, y_2, z_2)\in X'$, touches $C+{\bf x}_1$ at some relative interior points of $F_3+{\bf x}_1$. Clearly we have
$$y_2\ge y_1-2 \qquad {\rm and} \qquad z_2\ge z_1-2\ge y_1-2.$$
If $x_1>1.5$, by
$$\left\{
\begin{array}{ll}
(x_2-x_1)-(y_2-y_1)-(z_2-z_1)=4,&\\
x_1+y_1+z_1=4&
\end{array}\right.$$
we get
$$x_2-y_2-z_2=4+x_1-(y_1+z_1)=2x_1>3,$$
which contradicts the assumption ${\bf x}_2\in P_2$. Thus, in this case we have $x_1\le 1.5$ and therefore
$$z_1\ge 2-{{x_1}\over 2}\ge {5\over 4}. \eqno (4.40)$$

Now we consider two subcases.

\medskip\noindent
{\bf Subcase 1.1.} $y_2\ge z_1-2$. Let $S_7$ denote the halfspace $\{ (x,y,z):\ (x-x_1)-(y-y_1)-(z-z_1)\le 2\}$ and define
$$P_{10}=P_1\cap S_7.$$
By routine arguments based on Lemma 2.3 and Lemma 2.4, it can be deduced that
$$P_{10}\setminus \{ C+{\bf x}_1\}\subset D'_1(C,X').$$

We observe that $P_{10}$ is independent of ${\bf x}_2$, and
$${\rm vol}(P_{10}\setminus \{ C+{\bf x}_1\})\ge {\rm vol}(P_{10}\setminus \{ C+\overline{{\bf x}_1}\}),$$
where $\overline{{\bf x}_1}=(x_1, 2-{1\over 2}x_1, 2-{1\over 2}x_1)$. Let $T_{10}$ denote the tetrahedron with vertices $(1, 2-x_1, 0)$, $(x_1, 2-x_1, 0)$, $(1, 3-2x_1,0)$, and $(1, 2-x_1, 1-x_1)$, and let $T_{11}$
denote the tetrahedron with vertices $(1,0,2-x_1)$, $(x_1, 0,2-x_1)$, $(1, 1-x_1, 2-x_1)$ and $(1,0, 3-2x_1)$. Clearly, all $T_1$, $T_4$, $T_{10}$ and $T_{11}$ are congruent to each others. For convenience, we write ${\bf x}_2
=(2, 1-x_1, 1-x_1)$, $X^*=\{ {\bf o}, \overline{{\bf x}_1}, {\bf x}_2\}$,
$$R_5=\{ (x,y,z)\in \partial (C+{\bf x}_2):\ x\le 2,\ y\ge 1-x_1,\ z\ge 1-x_1\},$$ and define
$$P_{11}=\bigcup_{{\bf x}\in R_5}s(C, X^*, -{\bf e}_1, {\bf x}).$$

Then by Lemma 2.3 and Lemma 2.4 it can be verified that
$$P_{11} \subset D'_1(C,X')\cup T_{10}\cup T_{11}\setminus \{ T_1\cup T_4\}$$
and therefore
$${\rm vol}(D'_1(C, X'))\ge {\rm vol}(P_{11}).$$
In fact, the last equality holds only if $x_1=1$. By reflection, it can be seen that $P_{11}$ is certain $D'_1(C,X')$ estimated in Subcase 2.1 of Lemma 4.4. Thus, in this subcase we obtain
$${\rm vol}(D'_1(C,X'))> {1\over 6} \left({5\over 9}-{4\over 9}\sqrt{1\over {10}}\right).\eqno (4.41)$$

\medskip\noindent
{\bf Subcase 1.2.} $y_2\le z_1-2.$ Then the assumption ${\rm int}(F_3+{\bf x}_1)\cap J\not=\emptyset $ implies
$$y_1-2\le y_2\le z_1-2.$$
We define $T_0'$ to be the orthogonal tetrahedron with vertices $(1,y_2+1,z_1-1)$, $(2, y_2+1, z_1-1)$, $(1,y_2+2,z_1-1)$ and $(1, y_2+1, z_1-2)$, define $S_9$ to be the halfspace $\{ (x,y,z):\ z\ge z_2-1\}$, and define
$$P_3=P_1\cap T_0'\cap S_9.$$
It can be verified by Lemma 2.3 and Lemma 2.4 that
$$P_3\subseteq D'_1(C,X').$$

Similar to Subcase 2.2 of Lemma 4.4, by $x_2-y_2+z_2\le 4$, $(x_2-x_1)-(y_2-y_1)-(z_2-z_1)=4$ and $x_1+y_1+z_1=4$
we get
$$z_2\le 2-x_1\le 1.$$
Then, by the assumptions on ${\bf x}_1$ and ${\bf x}_2$ (in particular (4.38)) one can deduce
$$-{3\over 5}\le 3-2\beta \le 3-2z_1\le 1-z_1-y_2\le 3-y_1-z_1=x_1-1<{3\over 5},$$
\begin{eqnarray*}
{\rm vol}(T'_0\setminus \{S_2\cap S_9\})&\hspace{-0.2cm}\le &\hspace{-0.2cm}\max \left\{ {1\over 6}(2-z_1)^3,\  {1\over 6}(1-z_1+z_2)^3\right\}
\le {1\over 6}(2-z_1)^3,\\
{\rm vol}(T'_0\cap S_{10})&\hspace{-0.2cm}\le &\hspace{-0.2cm}{1\over 2}\left({{1-z_1-y_2}\over {\sqrt{2}}}\right)^2(z_1-1)\le {1\over 6}\times {3\over 2}\left( {3\over 5}\right)^2(z_1-1),\\
{\rm vol}(T'_0\cap S_{11})&\hspace{-0.2cm}\le &\hspace{-0.2cm}{1\over 6}\times 2\left({{z_1-y_2-2}\over 2}\right)^3\le {1\over 6}\times {1\over 4}(z_1-y_1)^3,
\end{eqnarray*}
and therefore
$${\rm vol}(P_3)\ge {1\over 6}\left(1-(1+y_2)^3-(2-z_1)^3-{3\over 2}\left( {3\over 5}\right)^2(z_1-1)-{1\over 4}(z_1-y_1)^3\right).$$

Thus, by (4.37), (4.40) and routine computations we get
\begin{eqnarray*}
\hspace{0.6cm}{\rm vol}(D'_1(C,X'))&\hspace{-0.2cm}\ge &\hspace{-0.2cm}{\rm vol}(T_1\cup T_2\cup T_3\cup T_4\cup P_3)\\
&\hspace{-0.2cm}\ge &\hspace{-0.2cm}{1\over 6}\left( 1+2(x_1-1)^3+(y_1-1)^3-(2-z_1)^3-{3\over 2}\left( {3\over 5}\right)^2(z_1-1)-{1\over 4}(z_1-y_1)^3\right)\\
&\hspace{-0.2cm}\ge &\hspace{-0.2cm}{1\over 6}\left( 1+2(x_1-1)^3-(2-z_1)^3-{3\over 2}\left( {3\over 5}\right)^2(z_1-1)-{1\over 4}(z_1-1)^3\right)\\
&\hspace{-0.2cm}\ge &\hspace{-0.2cm}{1\over 6}\left( 1-(2-\beta )^3-{3\over 2}\left( {3\over 5}\right)^2(\beta-1)-{1\over 4}(\beta -1)^3\right)\\
&\hspace{-0.2cm}> &\hspace{-0.2cm}{1\over 6}\left({5\over 9}-{4\over 9}\sqrt{1\over {10}}\right).\hspace{9.3cm} (4.42)
\end{eqnarray*}

As a conclusion of (4.41) and (4.42), in this case we obtain
$${\rm vol}(D'_1(C,X'))> {1\over 6} \left({5\over 9}-{4\over 9}\sqrt{1\over {10}}\right).$$

\medskip\noindent
{\bf Case 2.} ${\rm int}(F_0)\cap J\not= \emptyset .$ Assume that $C+{\bf x}_2$, where ${\bf x}_2=(2, y_2, z_2)\in
X'$, touches $C$ at some interior points of $F_0$. Similar to Case 3 of Lemma 4.4, we divide the region $\{(y_2, z_2):\ -1\le y_2\le 2,\ -1\le z_2\le 2\}$ into ten parts as illustrated in Figure 7 and consider the corresponding subcases.

\begin{figure}[hc]
\includegraphics[height=5.8cm,width=5.5cm,angle=0]{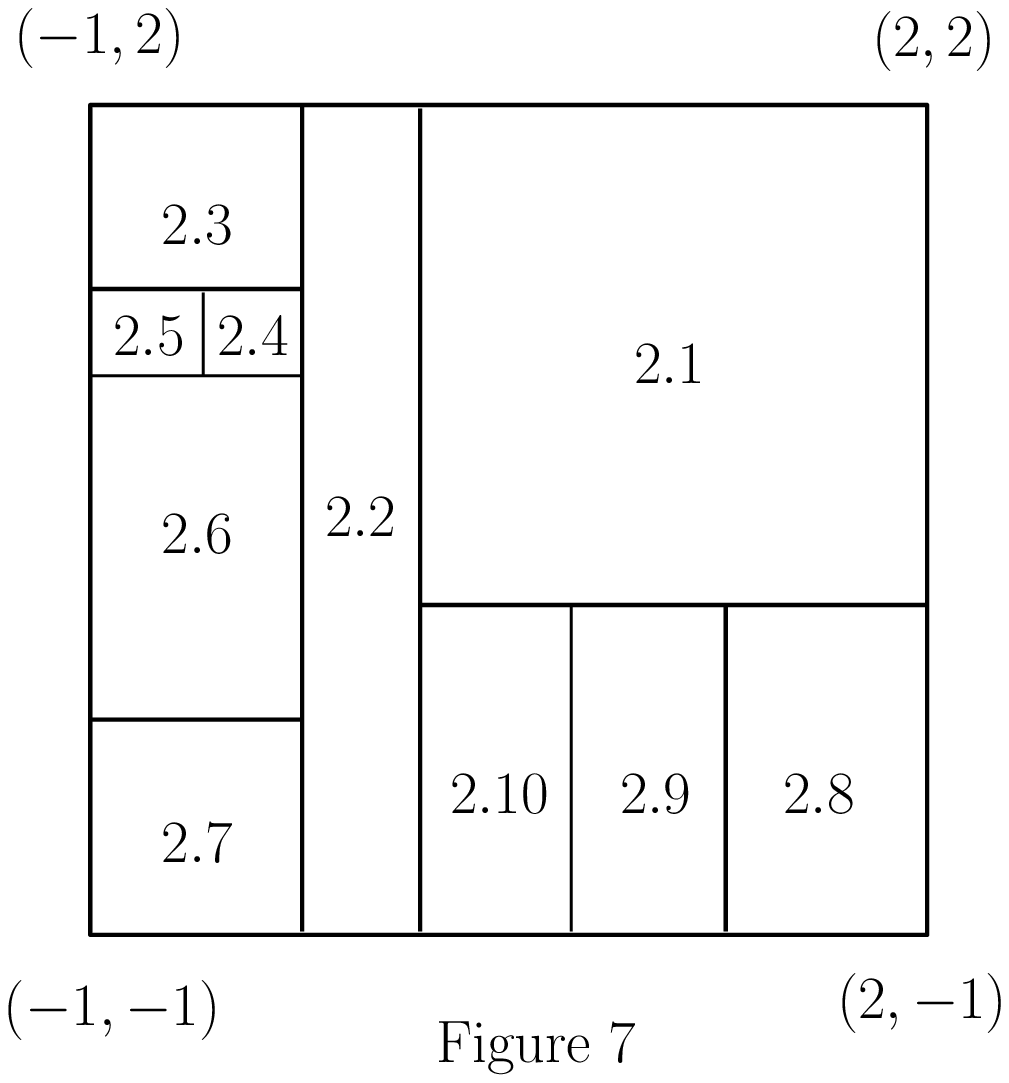}
\end{figure}

\medskip\noindent
{\bf Subcase 2.1.} $y_2\ge z_1-2$ {\it and} $z_2\ge z_1-2$. Then, it can be deduced that the hyperplane $(x-x_1)-(y-y_1)-(z-z_1)=2$ separates $C+{\bf x}_1$ and $C+{\bf x}_2$. Therefore, similar to Case 1, we have
$x_1\le {3\over 2}.$ Let $S_7$ denote the halfspace $\{ (x,y,z):\ (x-x_1)-(y-y_1)-(z-z_1)\le 2\}$ and define
$$P_{10}=P_1\cap S_7.$$
By routine arguments based on Lemma 2.3 and Lemma 2.4, it can be deduced that
$$P_{10}\setminus \{ C+{\bf x}_1\}\subseteq D'_1(C,X').$$

It is easy to see that $P_{10}$ is independent of ${\bf x}_2$, up to ${\bf x}_2$ satisfying $x_2=2$, $y_2\ge z_1-2$
and $z_2\ge z_1-2$, and
$${\rm vol}(P_{10}\setminus \{ C+{\bf x}_1\})\ge {\rm vol}(P_{10}\setminus \{C+\overline{{\bf x}_1}\}),$$
where $\overline{{\bf x}_1}=(x_1, 2-{1\over 2}x_1, 2-{1\over 2}x_1).$ Then, by arguments similar to Subcase 1.1,
it can be proved that
$${\rm vol}(D'_1(C, X')) > {1\over 6}\left({5\over 9}-{4\over 9}\sqrt{1\over {10}}\right).\eqno (4.43)$$

\medskip\noindent
{\bf Subcase 2.2.} $y_1-2\le y_2\le z_1-2$. Let $P_3$ be the polytope defined in Subcase 1.2, here we also have
$$P_3\subseteq D'_1(C,X').$$
If $z_1-z_2\le 2$, then we get
$$(2-x_1)-(y_2+1-y_1)-(z_2+1-z_1)\ge 2,$$
which implies
$$z_1-z_2\ge 2+(y_2-y_1)+x_1\ge x_1\ge 1.$$
Thus we have
$${\rm vol}(P_3)\ge {1\over 6}\left(1-(1+y_2)^3-(2-z_1)^3-{3\over 2}\left({3\over 5}\right)^2(z_1-1)-{1\over 4}(z_1-y_1)^3\right).$$

On the other hand, by the assumption ${\bf x}_1\in \triangle_3\cap H$ and (4.38) we get
$$1.2075\le z_1\le \beta .$$
Thus, by (4.37) and routine computations we get
\begin{eqnarray*}
\hspace{0.6cm}{\rm vol}(D'_1(C,X'))&\hspace{-0.2cm}\ge &\hspace{-0.2cm}{\rm vol}(T_1\cup T_2\cup T_3\cup T_4\cup P_3)\\
&\hspace{-0.2cm}\ge &\hspace{-0.2cm}{1\over 6}\left( 1+2(x_1-1)^3+(y_1-1)^3-(2-z_1)^3-{3\over 2}\left({3\over 5}\right)^2(z_1-1)-{1\over 4}(z_1-y_1)^3\right)\\
&\hspace{-0.2cm}\ge &\hspace{-0.2cm}{1\over 6}\left( 1+2(x_1-1)^3-(2-z_1)^3-{3\over 2}\left({3\over 5}\right)^2(z_1-1)-{1\over 4}(z_1-1)^3\right)\\
&\hspace{-0.2cm}\ge &\hspace{-0.2cm}{1\over 6} \left( 1-(2-\beta )^3-{3\over 2}\left({3\over 5}\right)^2(\beta -1)-{1\over 4}(\beta -1)^3\right)\\
&\hspace{-0.2cm}> &\hspace{-0.2cm} {1\over 6}\left({5\over 9}-{4\over 9}\sqrt{1\over {10}}\right).\hspace{9.3cm}(4.44)
\end{eqnarray*}

\medskip\noindent
{\bf Subcase 2.3.} {\it $-1\le y_2\le y_1-2$ and $z_2\ge 1.16$}. Then, by ${\bf x}_1\in \triangle_3\cap H$ and
$$2-y_2+z_2\le 4$$
one can deduce
$$1.16\le z_2\le y_2+2\le y_1,$$
$$-0.84\le y_2\le y_1-2\le -0.5,$$
and
$$-0.68\le y_2+z_2-1\le 2y_1-3\le 0.$$
Let $T_0^*$ denote the orthogonal tetrahedron with vertices $(1,y_2+1, z_2-1)$, $(2, y_2+1, z_2-1)$, $(1,y_2,z_2-1)$ and $(1,y_2+1,z_2)$, and let $P_4$ be the polytope defined in Subcase 3.3 of Lemma 4.4. Since $x_1\ge 1$ holds in this case, the $S_8$ is no longer effective in the definition of $P_4$. It can be verified by Lemma 2.3 and Lemma 2.4 that
$$P_4\subset D'_1(C,X').$$

Recall that $S_1=\{ (x,y,z):\ y\ge 0\}$, $S'_6=\{ (x,y,z):\ x+y-z\ge 2\}$ and $S_{10}=\{ (x,y,z):\ x-y-z\ge 1\}$, we have
$${\rm vol}(T_0^*\cap S'_6)=2\times {1\over 6}\left({{2+y_2-z_2}\over 2}\right)^3={1\over 6}\times {1\over 4}\ (2+y_2-z_2)^3,$$
$${\rm vol}(T_0^*\cap S_1\cap S_{10})\le {1\over 6}\times {3\over 2}\ (y_2+z_2-1)^2(y_2+1)$$
and therefore, based on the assumptions on ${\bf x}_1$ and ${\bf x}_2$,
\begin{eqnarray*}
{\rm vol}(P_4) &\hspace{-0.2cm}\ge &\hspace{-0.2cm} {1\over 6}\left(1-(z_2-1)^3-{1\over 4}(2+y_2-z_2)^3+y_2^3-{3\over 2}(y_2+z_2-1)^2(y_2+1)\right)\\
&\hspace{-0.2cm}\ge &\hspace{-0.2cm}{1\over 6} \left(1-(y_2+1)^3+y_2^3-{3\over 2}(y_2+z_2-1)^2(y_2+1)\right)\\
&\hspace{-0.2cm}\ge &\hspace{-0.2cm}{1\over 6} \left(1-(y_2+1)^3+y_2^3-{3\over 2}\times 0.68^2(y_2+1)\right)\\
&\hspace{-0.2cm}\ge &\hspace{-0.2cm}{1\over 6} \left(1-\left( 1-0.84 \right)^3-0.84^3-{3\over 2}\times 0.68 ^2(1-0.84)\right).
\end{eqnarray*}

Then, together with (4.37), we get
\begin{eqnarray*}
\hspace{1cm}{\rm vol}(D'_1(C,X'))&\hspace{-0.2cm}\ge &\hspace{-0.2cm}{\rm vol}(T_1\cup T_2\cup T_3\cup T_4\cup P_4)\\
&\hspace{-0.2cm}\ge &\hspace{-0.2cm} {1\over 6}\left(1+{{4\sqrt{2}+2}\over {25+22\sqrt{2}}}-\left( 1-0.84 \right)^3-0.84^3-{3\over 2}\times 0.68 ^2(1-0.84)\right)\\
&\hspace{-0.2cm}> &\hspace{-0.2cm}{1\over 6}\left({5\over 9}-{4\over 9}\sqrt{1\over {10}}\right). \hspace{8.9cm}(4.45)
\end{eqnarray*}

\medskip\noindent
{\bf Subcase 2.4.} {\it $-0.8\le y_2\le y_1-2$ and $1\le z_2\le 1.16$}. These conditions and ${\bf x}_1\in \triangle_3\cap H$ imply
$$-0.8\le y_2\le -0.5.$$
Let $P_4$ be the polytope defined in the previous subcase. Then we have
$$P_4\subset D'_1(C,X')$$
and
\begin{eqnarray*}
{\rm vol}(P_4) &\hspace{-0.2cm}\ge &\hspace{-0.2cm} {1\over 6}\left(1-(z_2-1)^3-{1\over 4}(2+y_2-z_2)^3+y_2^3-{3\over 2}(y_2+z_2-1)^2(y_2+1)\right)\\
&\hspace{-0.2cm}\ge &\hspace{-0.2cm}{1\over 6} \left(1-(1.16-1)^3-{1\over 4}(1+y_2)^3+y_2^3-{3\over 2}y_2^2(1+y_2)\right)\\
&\hspace{-0.2cm}\ge &\hspace{-0.2cm}{1\over 6} \left(1-0.16^3-{1\over 4}(1-0.8)^3-0.8^3-{3\over 2}\times 0.8^2(1-0.8)\right).
\end{eqnarray*}

Then, combined with (4.37), we get
\begin{eqnarray*}
\hspace{1.1cm}{\rm vol}(D'_1(C,X'))&\hspace{-0.2cm}\ge &\hspace{-0.2cm}{\rm vol}(T_1\cup T_2\cup T_3\cup T_4\cup P_4)\\
&\hspace{-0.2cm}\ge &\hspace{-0.2cm} {1\over 6}\left(1+{{4\sqrt{2}+2}\over {25+22\sqrt{2}}}-0.16^3-{1\over 4}(1-0.8)^3-0.8^3-{3\over 2}\times 0.8^2(1-0.8)\right)\\
&\hspace{-0.2cm}> &\hspace{-0.2cm}{1\over 6}\left({5\over 9}-{4\over 9}\sqrt{1\over {10}}\right). \hspace{8.9cm}(4.46)
\end{eqnarray*}

\medskip
\noindent
{\bf Subcase 2.5.} {\it $-1\le y_2\le \min\{ -0.8,\  y_1-2\}$ and $1\le z_2\le 1.16$}. Then, by the fact that $(2,y_2,z_2)\in P_2$ one can deduce
$$y_2+1\ge z_2-1.\eqno (4.47)$$
Let $P_3$ and $P_4$ be the polytopes defined in Subcases 1.2 and 2.3, respectively. It can be verified that
$${\rm int}(P_3)\cap {\rm int}(P_4)=\emptyset $$
and, by Lemma 2.3 and Lemma 2.4,
$$P_i\subset D'_1(C,X'),\qquad i=3,\ 4.$$

Then, it can be deduced that
$${\rm vol}(P_3)\ge {1\over 6}\left(1-(y_2+1)^3-{1\over 4}(z_1-y_2-2)^3-(1-z_1+z_2)^3-{3\over 2} (1-y_2-z_1)^2(z_1-z_2)\right).$$
In particular, when $1.35\le z_1\le \beta $, $-1\le y_2\le -0.8$ and $1\le z_2\le 1.16$, by routine computations we get
\begin{eqnarray*}
{\rm vol}(P_3)&\hspace{-0.2cm}\ge &\hspace{-0.2cm} {1\over 6}\left(1-0.2^3-{1\over 4}(z_1-1)^3-(1-z_1+z_2)^3-{3\over 2} (2-z_1)^2(z_1-z_2)\right)\\
&\hspace{-0.2cm}\ge &\hspace{-0.2cm} {1\over 6}\left(1-0.2^3-{1\over 4}(z_1-1)^3-(2.16-z_1)^3-{3\over 2} (2-z_1)^2(z_1-1.16)\right)\\
&\hspace{-0.2cm}\ge &\hspace{-0.2cm} {1\over 6}\left(1-0.2^3-{1\over 4}(z_1-1)^3-(2.16-z_1)^3-{3\over 2}\times 0.65^2(z_1-1.16)\right)\\
&\hspace{-0.2cm}>&\hspace{-0.2cm}{1\over 6}\times {3\over {10}}
\end{eqnarray*}
and therefore
$${\rm vol}(D'_1(C, X'))\ge {\rm vol}(T_1\cup T_2\cup T_3\cup T_4\cup P_3)\ge {1\over 6}\left({{4\sqrt{2}+2}\over {25+22\sqrt{2}}}+{3\over {10}}\right)> {1\over 6}\left( {5\over 9}-{4\over 9}\sqrt{1\over {10}}\right).\eqno (4.48)$$

\begin{figure}[hc]
\includegraphics[height=5.2cm,width=11.5cm,angle=0]{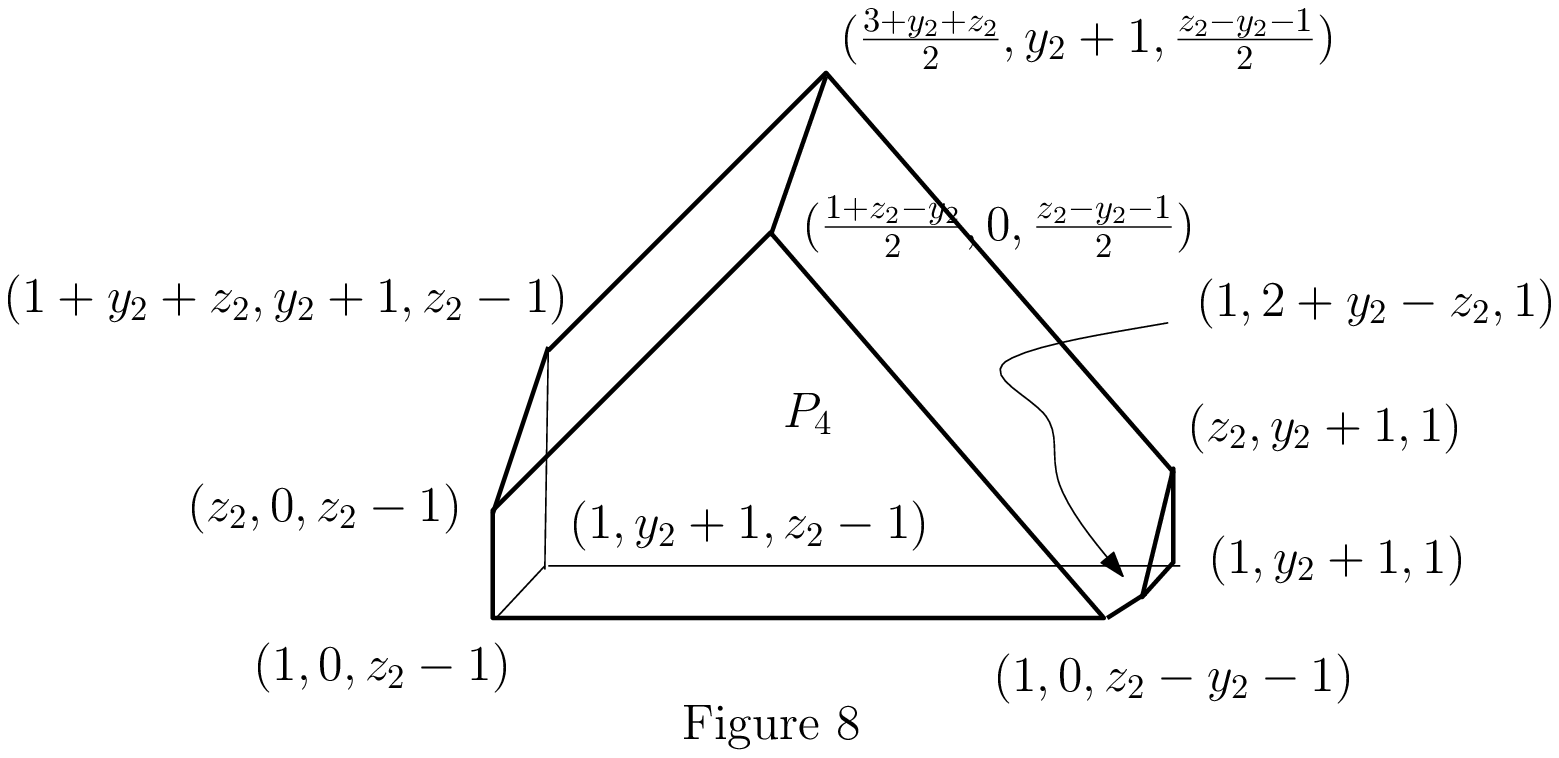}
\end{figure}

In fact, as shown in Figure 8, $P_4$ is the polytope with ten vertices $(1, y_2+1, z_2-1)$, $(1+y_2+z_2, y_2+1, z_2-1)$, $({{3+y_2+z_2}\over 2}, y_2+1, {{z_2-y_2-1}\over 2})$, $(z_2, y_2+1, 1)$, $(1, y_2+1, 1),$ $(1, 0, z_2-1)$, $(z_2, 0, z_2-1)$, $({{1+z_2-y_2}\over 2}, 0, {{z_2-y_2-1}\over 2})$, $(1, 0, z_2-y_2-1)$ and $(1, 2+y_2-z_2, 1)$. We observe that $P_4$ can be obtained by cutting of two tetrahedra and one prismoid from the tetrahedron with vertices $({{3+y_2+z_2}\over 2}, y_2+1, {{z_2-y_2-1}\over 2})$, $(1, y_2+1, -y_2-1)$, $(1, y_2+1, z_2)$ and $(1, {{y_2-z_2+1}\over 2}, {{z_2-y_2-1}\over 2})$, where the last tetrahedron is the union of two tetrahedra, both of them are homothetic to $T_0$ with ratio ${{y_2+z_2+1}\over 2}$. Then, by (4.47), $-1\le y_2\le -0.8$ and $1\le z_2\le 1.16$, we get
\begin{eqnarray*}
\hspace{1.3cm}{\rm vol}(P_4)&\hspace{-0.2cm} =&\hspace{-0.2cm} {1\over 6} \left( {1\over 4}(z_2+(y_2+1))^3-{1\over 4}(z_2-(y_2+1))^3-(z_2-1)^3 -(y_2+1) \left((y_2+z_2)^2 \right.\right.\\
&\hspace{-0.2cm}&\hspace{-0.2cm} \left.+(z_2-1)^2+(y_2+z_2)(z_2-1)\right)\biggr)\\
&\hspace{-0.2cm} \ge &\hspace{-0.2cm} {1\over 6} \left( {3\over 2}z^2_2(y_2+1)+{1\over 2}(y_2+1)^3-(y_2+1)(z_2-1)^2 -(y_2+1) \left((y_2+z_2)^2
\right.\right.\\
&\hspace{-0.2cm}&\hspace{-0.2cm} \left.+(z_2-1)^2+(y_2+z_2)(z_2-1)\right)\biggr)\\
&\hspace{-0.2cm}\ge &\hspace{-0.2cm}{{y_2+1}\over 6}\left({3\over 2}z_2^2+{1\over 2}(1-0.8)^2-2(z_2-1)^2-2(z_2-0.8)^2\right)\\
&\hspace{-0.2cm}\ge &\hspace{-0.2cm}{1\over 6}\times 1.4\times (y_2+1).\hspace{9.7cm}(4.49)
\end{eqnarray*}

Similarly, when $1.2075\le z_1\le 1.35$, noticing $y_2+1+z_1=(y_2+1+z_2)+(z_1-z_2)$ and $y_2+1-z_1+2z_2=
(y_2+1+z_2)-(z_1-z_2)$, we get
\begin{eqnarray*}
\hspace{1.3cm}{\rm vol}(P_3)&\hspace{-0.2cm} \ge&\hspace{-0.2cm} {1\over 6} \left( {1\over 4}(y_2+1+z_1)^3-{1\over 4}(y_2+1-z_1+2z_2)^3-3(z_1-z_2)\left((y_2+z_1)^2+(y_2+1)^2\right)\right)\\
&\hspace{-0.2cm}\ge &\hspace{-0.2cm}{{z_1-z_2}\over 6}\left({3\over 2}(y_2+1+z_2)^2+{1\over 2}(z_1-z_2)^2-3(y_2+z_1)^2-3(y_2+1)^2\right)\\
&\hspace{-0.2cm}\ge &\hspace{-0.2cm}{1\over 6}\left({3\over 2}-{5\over 2}(z_1-1)^2\right)(z_1-z_2).\hspace{8cm}(4.50)
\end{eqnarray*}
At the same time, we have
$${3\over 2}-{5\over 2}(z_1-1)^2\le 1.4.\eqno (4.51)$$

Then, by (4.37), (4.49)-(4.51) and (4.47) we get
\begin{eqnarray*}
\hspace{1.5cm}{\rm vol}(D'_1(C,X'))&\hspace{-0.2cm}\ge &\hspace{-0.2cm}{\rm vol}(T_1\cup T_2\cup T_3\cup T_4\cup P_3\cup P_4)\\
&\hspace{-0.2cm}\ge &\hspace{-0.2cm} {1\over 6}\left( {{4\sqrt{2}+2}\over {25+22\sqrt{2}}}+\left({3\over 2}-{5\over 2}(z_1-1)^2\right)
(z_1-z_2)+1.4 (y_2+1)\right)\\
&\hspace{-0.2cm}\ge &\hspace{-0.2cm} {1\over 6}\left( {{4\sqrt{2}+2}\over {25+22\sqrt{2}}}+\left({3\over 2}-{5\over 2}(z_1-1)^2\right)
(z_1-z_2+y_2+1)\right)\\
&\hspace{-0.2cm}\ge &\hspace{-0.2cm} {1\over 6}\left({{4\sqrt{2}+2}\over {25+22\sqrt{2}}}+\left({3\over 2}-{5\over 2}(z_1-1)^2\right)
(z_1-1)\right)\\
&\hspace{-0.2cm}> &\hspace{-0.2cm}{1\over 6}\left({5\over 9}-{4\over 9}\sqrt{1\over {10}}\right).\hspace{8.4cm}(4.52)
\end{eqnarray*}

\medskip\noindent
{\bf Subcase 2.6.} {\it $-1\le y_2\le y_1-2$ and $y_1-2\le z_2\le 1$}. Then we get
$$P_3\subset D'_1(C,X'),$$
where $P_3$ was defined in Subcase 1.2, and
$${\rm vol}(P_3)\ge {1\over 6}\left(1-(1+y_2)^3-{1\over 4}(z_1-y_2-2)^3-(2-z_1)^3-
{3\over 2}(1-y_2-z_1)^2(z_1-1)\right).$$

By routine computations one can deduce
$$(1+y_2)^3+{1\over 4} (z_1-y_2-2)^3\le  (z_1-1)^3,$$
$$x_1-1=3-y_1-z_1\le 1-y_2-z_1\le 2-z_1$$
and
$$(1-y_2-z_1)^2(z_1-1)\le (2-z_1)^2(z_1-1).$$
On the other hand, when ${\bf x}_1\in \triangle_3\cap H$, by (4.38) we get
$$z_1\ge 2-{{x_1}\over 2}\ge 1.2075.$$
Thus, by the assumptions on ${\bf x}_1$, ${\bf x}_2$ and (4.37) we get
\begin{eqnarray*}
\hspace{0.8cm}{\rm vol}(D'_1(C,X'))&\hspace{-0.2cm}\ge &\hspace{-0.2cm} {\rm vol}(T_1\cup T_2\cup T_3\cup T_4\cup P_3)\\
&\hspace{-0.2cm}\ge &\hspace{-0.2cm}{1\over 6}\left(1+2(x_1-1)^3+(y_1-1)^3+(z_1-1)^3-(1+y_2)^3-{1\over 4}(z_1-y_2-2)^3\right.\\
&\hspace{-0.2cm}&\hspace{-0.2cm}-\left. (2-z_1)^3-{3\over 2}(1-y_2-z_1)^2(z_1-1)\right)\\
&\hspace{-0.2cm}\ge &\hspace{-0.2cm} {1\over 6}\left( 1+2(x_1-1)^3+(y_1-1)^3-(2-z_1)^3-{3\over 2}(2-z_1)^2(z_1-1)\right)\\
&\hspace{-0.2cm}= &\hspace{-0.2cm} {1\over 6}\left( 1+2(x_1-1)^3+(y_1-1)^3-{1\over 2}(2-z_1)^2(1+z_1)\right)\\
&\hspace{-0.2cm}\ge &\hspace{-0.2cm} {1\over 6}\left( 1+{{4\sqrt{2}+2}\over {25+22\sqrt{2}}}-(z_1-1)^3-{1\over 2}(2-z_1)^2(1+z_1)\right)\\
&\hspace{-0.2cm}\ge &\hspace{-0.2cm} {1\over 6}\left( 1+{{4\sqrt{2}+2}\over {25+22\sqrt{2}}}-(1.2075-1)^3-{1\over 2}(2-1.2075)^2(1+1.2075)\right)\\
&\hspace{-0.2cm}> &\hspace{-0.2cm} {1\over 6}\left({5\over 9}-{4\over 9}\sqrt{1\over {10}}\right). \hspace{9.2cm}(4.53)
\end{eqnarray*}

\medskip\noindent
{\bf Subcase 2.7.} {\it $-1\le y_2\le y_1-2$ and $z_2<y_1-2$}. First, by the assumption ${\bf x}_2=(2,y_2,z_2)\in P_2$ one can deduce
$$2-y_2-z_2<3$$
and therefore
$$y_2+z_2> -1.\eqno (4.54)$$
On the other hand, by the assumption $y_2\le y_1-2$, $z_2<y_1-2$ and ${\bf x}_1\in \triangle_3\cap H$ we get
$$y_2+z_2<2y_1-4\le -1,$$
which contradicts (4.54). Therefore this subcase can't happen.

\medskip\noindent
{\bf Subcase 2.8.} {\it $-1\le z_2\le z_1-2$ and $y_2\ge 1.2$}. By the assumption ${\bf x}_2=(2,y_2,z_2)\in P_2$ we get
$$2+y_2-z_2< 4.$$
Thus one can deduce
$$1.2\le y_2< z_2+2\le z_1,$$
which implies
$$-0.8\le z_2\le z_1-2$$
and
$$z_2+0.2\le y_2+z_2-1\le 2z_2+1\le 2z_1-3\le 2\beta -3.$$
Let $T_0^\bullet $ denote the orthogonal tetrahedron with vertices $(1,y_2-1, z_2+1)$, $(2, y_2-1, z_2+1)$, $(1,y_2,z_2+1)$ and $(1,y_2-1,z_2)$, and let $P_6$ be the polytope defined in Subcase 3.7 of Lemma 4.4.
Now, that is
$$P_6=P_1\cap T_0^\bullet .$$
By Lemma 2.3 and Lemma 2.4 it can be verified that
$$P_6\subset D'_1(C,X').$$

We recall $S_2=\{ (x,y,z):\ z\ge 0\}$, $S_{10}=\{ (x,y,z):\  x-y-z\ge 1\}$ and $S_{11}=\{ (x,y,z):\ x-y+z\ge 2\}$. By routine computations we get
$${\rm vol}(T^\bullet_0\cap S_{11})={1\over 6}\times {1\over 4}(2-y_2+z_2)^3,$$
$${\rm vol}(T^\bullet_0\cap S_2\cap S_{10})\le {1\over 6}\times {3\over 2}(y_2+z_2-1)^2(z_2+1),$$
$${\rm vol}(P_6) \ge {1\over 6}\left(1-(y_2-1)^3-{1\over 4}(2-y_2+z_2)^3+z_2^3-{3\over 2}(y_2+z_2-1)^2(z_2+1)\right)$$
and therefore
\begin{eqnarray*}
{\rm vol}(D'_1(C,X'))&\hspace{-0.2cm}\ge &\hspace{-0.2cm}{\rm vol}(T_1\cup T_2\cup T_3\cup T_4\cup P_6)\\
&\hspace{-0.2cm}\ge &\hspace{-0.2cm}{1\over 6}\left(1+2(x_1-1)^3+(y_1-1)^3+(z_1-1)^3-(y_2-1)^3-{1\over 4}(2-y_2+z_2)^3+z_2^3\right.\\
&\hspace{-0.2cm}&\hspace{-0.2cm}\left. -{3\over 2}(y_2+z_2-1)^2(z_2+1)\right)\\
&\hspace{-0.2cm}\ge &\hspace{-0.2cm} {1\over 6} \left(1+2(x_1-1)^3+(y_1-1)^3+(z_1-1)^3-(z_2+1)^3+z_2^3-{3\over 2}(y_2+z_2-1)^2(z_2+1)\right).
\end{eqnarray*}

When $-0.8\le z_2\le -0.4$, we have
$$z_2+0.2\le y_2+z_2-1\le 2z_2+1\le -(z_2+0.2)$$
and
\begin{eqnarray*}
\hspace{0.5cm}{\rm vol}(D'_1(C, X'))&\hspace{-0.2cm}\ge &\hspace{-0.2cm} {1\over 6} \left(1+2(x_1-1)^3+(y_1-1)^3+(z_1-1)^3-(z_2+1)^3+z_2^3-{3\over 2}(z_2+0.2)^2(z_2+1)\right)\\
&\hspace{-0.2cm}\ge &\hspace{-0.2cm} {1\over 6}\left( 1+{{4\sqrt{2}+2}\over {25+22\sqrt{2}}}-\left({1\over 5}\right)^3-\left({4\over 5}\right)^3
-{3\over 2}\left({3\over 5}\right)^2{1\over 5}\right)\\
&\hspace{-0.2cm}> &\hspace{-0.2cm}{1\over 6}\times {1\over 2}. \hspace{11.2cm}(4.55)
\end{eqnarray*}
When $-0.4\le z_2\le z_1-2$, we get
$$-(2z_2+1)\le z_2+0.2\le y_2+z_2-1\le 2z_2+1$$
and
\begin{eqnarray*}
\hspace{0.4cm}{\rm vol}(D'_1(C,X'))&\hspace{-0.2cm}\ge &\hspace{-0.2cm} {1\over 6} \left( 1+2(x_1-1)^3+(y_1-1)^3+(z_1-1)^3-(z_2+1)^3+z_2^3-{3\over 2}(2z_2+1)^2(z_2+1)\right)\\
&\hspace{-0.2cm}\ge &\hspace{-0.2cm} {1\over 6} \left(1+2(x_1-1)^3+(y_1-1)^3-(2-z_1)^3-{3\over 2}(2z_1-3)^2(z_1-1)\right)\\
&\hspace{-0.2cm}\ge &\hspace{-0.2cm} {1\over 6} \left(1+2(x_1-1)^3-(2-z_1)^3-{3\over 2}(2z_1-3)^2(z_1-1)\right)\\
&\hspace{-0.2cm}\ge &\hspace{-0.2cm} {1\over 6} \left(1-(2-\beta )^3-{3\over 2}(2\beta -3)^2(\beta -1)\right)\\
&\hspace{-0.2cm}> &\hspace{-0.2cm}{1\over 6}\times {1\over 2}. \hspace{11.3cm}(4.56)
\end{eqnarray*}

As a conclusion of (4.55) and (4.56), in this subcase we have
$${\rm vol}(D'_1(C, X'))> {1\over 6}\times {1\over 2}.$$

\medskip\noindent
{\bf Subcase 2.9.} {\it $-1\le z_2\le z_1-2$ and $1\le y_2\le 1.2$}. Let $P_6$ be the polytope defined in the
previous subcase and define
\begin{eqnarray*}
P_7&\hspace{-0.2cm}=&\hspace{-0.2cm}P_1\cap \left\{ (x,y,z):\ 1\le x\le 2,\ y_2-1 \le y\le 1, \ z_2+1\le z\le z_1-1,\right.\\
&\hspace{-0.2cm}&\hspace{-0.2cm}\left. (x-1)-(y-y_1+1)+(z-z_2-1)\le 1\right\}.
\end{eqnarray*}
Clearly, we have
$${\rm int}(P_6)\cap {\rm int}(P_7) = \emptyset .$$
By Lemmas 2.3 and 2.4 it can be shown that
$$P_i\subseteq D'_1(C, X')$$
holds for both $i=6$ and $7$.

For convenience, we write $$\gamma (z_1)=\min\{1,\ z_1-0.5,\ 2.5-z_1\}.$$
It can be verified that $\gamma (z_1)\ge 0.7$ and $P_7$ contains $(1,1, z_2+1)$, $(1, y_2-1, z_2+1)$,
$(1, y_2-1, z_1-1)$, $(1,1, z_1-1)$, $(1.5, 0.5, z_2+1)$ and $(1+\gamma (z_1), 0.5, z_1-1)$, and therefore
\begin{eqnarray*}
{\rm vol}(P_7)&\hspace{-0.2cm}\ge &\hspace{-0.2cm}{1\over 3} (2-y_2)(z_1-z_2-2)\gamma (z_1)+{1\over 6}(2-y_2)(z_1-z_2-2)(1.5-1)\\
&\hspace{-0.2cm}\ge &\hspace{-0.2cm}{1\over 6}\times 1.52 (z_1-z_2-2).
\end{eqnarray*}
In addition, by $2+y_2-z_2\le 4$ we get
\begin{eqnarray*}
{\rm vol}(P_6) &\hspace{-0.2cm}\ge &\hspace{-0.2cm} {1\over 6}\left(1-(y_2-1)^3-{1\over 4}(2-y_2+z_2)^3+z_2^3-{3\over 2}(y_2+z_2-1)^2(z_2+1)\right)\\
&\hspace{-0.2cm}\ge &\hspace{-0.2cm} {1\over 6}\left(1-(z_2+1)^3+z_2^3-{3\over 2}z_2^2(z_2+1)\right).
\end{eqnarray*}
Then we have
\begin{eqnarray*}
\hspace{0.7cm}{\rm vol}(D'_1(C, X'))&\hspace{-0.2cm}\ge &\hspace{-0.2cm}{\rm vol}(T_1\cup T_2\cup T_3\cup T_4\cup P_6\cup P_7)\\
&\hspace{-0.2cm}\ge &\hspace{-0.2cm}{1\over 6}\Bigl(1+2(x_1-1)^3+(y_1-1)^3+(z_1-1)^3+1.52(z_1-z_2-2)-(z_2+1)^3+z_2^3\\
&\hspace{-0.2cm}&\hspace{-0.2cm} -{3\over 2}z_2^2(z_2+1)\Bigr)\\
&\hspace{-0.2cm}\ge &\hspace{-0.2cm}{1\over 6}\left(1+2(x_1-1)^3+(y_1-1)^3-(2-z_1)^3-{3\over 2}(2-z_1)^2(z_1-1)
\right)\\
&\hspace{-0.2cm}\ge &\hspace{-0.2cm}{1\over 6}\left(1+{{4\sqrt{2}+2}\over {25+22\sqrt{2}}}-(z_1-1)^3-(2-z_1)^3-{3\over 2}(2-z_1)^2(z_1-1)
\right)\\
&\hspace{-0.2cm}>&\hspace{-0.2cm}{1\over 6}\left({5\over 9}-{4\over 9}\sqrt{1\over {10}}\right). \hspace{9.3cm}(4.57)
\end{eqnarray*}

\medskip\noindent
{\bf Subcase 2.10.} {\it $-1\le z_2\le z_1-2$ and $z_1-2\le y_2\le 1$}. Then, by $y_2+z_2>-1$ we get $y_2>-0.5$. We define
\begin{eqnarray*}
P_7&\hspace{-0.2cm}=&\hspace{-0.2cm}P_1\cap \left\{ (x,y,z):\ 1\le x\le 2,\ 0\le y\le 1, \ z_2+1\le z\le z_1-1,\right.\\
&\hspace{-0.2cm}&\hspace{-0.2cm}\left. (x-1)-(y-y_1+1)+(z-z_2-1)\le 1\right\},\\
P_8&\hspace{-0.2cm}=&\hspace{-0.2cm}P_1\cap \left\{ (x,y,z):\ 1\le x\le 2,\ 0\le y\le 1,\ 0\le z\le z_2+1\right\}\setminus \{ C+{\bf x}_2\},\\
\overline{P_7}&\hspace{-0.2cm}=&\hspace{-0.2cm}P_7\cap \{ (x,y,z):\ z=z_2+1\}\qquad {\rm and}\qquad \overline{P_8}=P_8\cap \{ (x,y,z):\ z=z_2+1\}.
\end{eqnarray*}
In addition, we define
\begin{eqnarray*}
P_{12}&\hspace{-0.2cm}=&\hspace{-0.2cm}P_1\cap \left\{ (x,y,z):\ 1\le x\le 2,\ 0\le y\le y_1-1, \ z_1-1\le z\le 1,\right.\\
&\hspace{-0.2cm}&\hspace{-0.2cm}\left. (x-1)-(y-y_1+1)+(z-z_2-1)\le 1\right\}.
\end{eqnarray*}
Clearly, ${\rm int}(P_7)$, ${\rm int}(P_8)$ and ${\rm int}(P_{12})$ are pairwise disjoint.
By Lemmas 2.3 and 2.4 it can be shown that
$$P_i\subseteq D'_1(C, X')$$
holds for all $i=7$, $8$ and $12$.

When $z_2\ge y_1-2$, we have
$$P_7=P_1\cap \left\{ (x,y,z):\ 1\le x\le 2,\ 0 \le y\le 1, \ z_2+1\le z\le z_1-1\right\}$$
and
$$P_{12}=P_1\cap \left\{ (x,y,z):\ 1\le x\le 2,\ 0 \le y\le y_1-1, \ z_1-1\le z\le 1\right\},$$
where the first one implies
$$\overline{P_8}\subseteq \overline{P_7}.$$
Now we consider ${\rm vol}(P_7)$, ${\rm vol}(P_8)$ and ${\rm vol}(P_{12})$ as functions of $x_1$, $y_1$, $z_1$, $y_2$ and $z_2$. It can be shown that
$${{\partial {\rm vol}(P_7)}\over {\partial z_2}}=-s(\overline{P_7}),$$
$${{\partial {\rm vol}(P_8)}\over {\partial z_2}}\le s(\overline{P_8}),$$
$${{\partial ({\rm vol}(P_7)+{\rm vol}(P_8))}\over {\partial z_2}}\le -s(\overline{P_7})+s(\overline{P_8})\le 0$$
and therefore
\begin{eqnarray*}
{\rm vol}(P_7)+{\rm vol}(P_8)&\hspace{-0.2cm}\ge &\hspace{-0.2cm}({\rm vol}(P_7)+{\rm vol}(P_8))\bigg|_{z_2=z_1-2}\ge {\rm vol}(P_8)\bigg|_{z_2=z_1-2} \\
&\hspace{-0.2cm}\ge &\hspace{-0.2cm}\left\{
\begin{array}{ll}
{1\over 6} \left({1\over 2} (z_1-{1\over 2})^3-2(z_1-1)^3-{1\over 2}({3\over 2}-z_1)^3\right), & \mbox{ if $z_1\le {3\over 2}$,}\\
{1\over 3}\left(1-{{z_1}\over 3}\right)^3, & \mbox{ if $z_1\ge {3\over 2}$},
\end{array}
\right.\\
&\hspace{-0.2cm}= &\hspace{-0.2cm}\left\{
\begin{array}{ll}
{1\over 6} \left(3\left({1\over 2}\right)^2(z_1-1)-(z_1-1)^3\right), & \mbox{ if $z_1\le {3\over 2}$,}\\
{1\over 3}\left(1-{{z_1}\over 3}\right)^3, & \mbox{ if $z_1\ge {3\over 2}.$}
\end{array}
\right.
\end{eqnarray*}
For example, when $y_2={1\over 2}$ and $z_2\le -{1\over 2}$, the set $P_8$ can be illustrated by the following figure.
\begin{figure}[hc]
\includegraphics[height=4.5cm,width=10cm,angle=0]{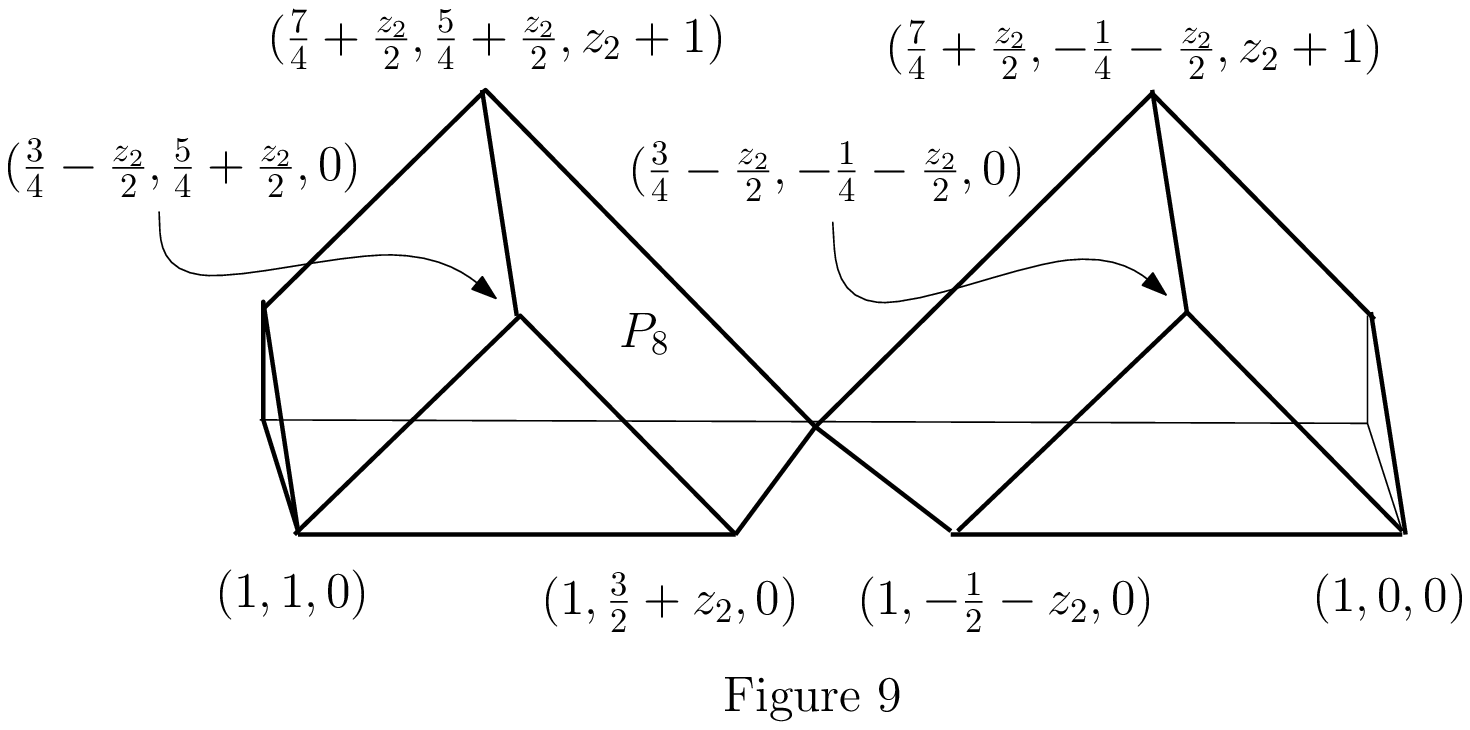}
\end{figure}

At the same time, we have
$${\rm vol}(P_{12})= \left\{
\begin{array}{ll}
\left(\left({1\over 2}\right)^3+{{z_1-0.5}\over 2}({3\over 2}-z_1)\right)(y_1-1)+{1\over 2} (y_1-1)^2(2-z_1), & \mbox{ if $z_1\le {3\over 2}$,}\\
{1\over 2}\left(2-z_1\right)^2(y_1-1)+{1\over 2}(y_1-1)^2(2-z_1), & \mbox{ if $z_1\ge {3\over 2}$}.
\end{array}
\right.$$

When $z_1\le {3\over 2}$, together with (4.37), we get
\begin{eqnarray*}
\hspace{2.4cm}{\rm vol}(D'_1(C, X'))&\hspace{-0.2cm}\ge &\hspace{-0.2cm}{\rm vol}(T_1\cup T_2\cup T_3\cup T_4\cup P_7\cup P_8\cup P_{12})\\
&\hspace{-0.2cm}\ge &\hspace{-0.2cm}{1\over 6}\left(2(x_1-1)^3+(y_1-1)^3+{3\over 4}(y_1-1)+{3\over 4} (z_1-1)\right)\\
&\hspace{-0.2cm}= &\hspace{-0.2cm}{1\over 6}\left(2(x_1-1)^3+(y_1-1)^3+{3\over 4}(2-x_1)\right)\\
&\hspace{-0.2cm}>&\hspace{-0.2cm}{1\over 6}\times {1\over 2}. \hspace{9.3cm}(4.58)
\end{eqnarray*}
When $z_1\ge {3\over 2}$, combined with (4.37), we get
\begin{eqnarray*}
\hspace{2.3cm}{\rm vol}(D'_1(C, X'))&\hspace{-0.2cm}\ge &\hspace{-0.2cm}{\rm vol}(T_1\cup T_2\cup T_3\cup T_4\cup P_7\cup P_8\cup P_{12})\\
&\hspace{-0.2cm}\ge &\hspace{-0.2cm}{1\over 6}\left(2(x_1-1)^3+(y_1-1)^3+(z_1-1)^3+2\left(1-{{z_1}\over 3}\right)^3\right.\\
&\hspace{-0.2cm}&\hspace{-0.2cm}+3(2-z_1)^2(y_1-1)+ 3(y_1-1)^2(2-z_1)\biggr)\\
&\hspace{-0.2cm}\ge &\hspace{-0.2cm}{1\over 6}\left(2(x_1-1)^3+3(2-z_1)^2(y_1-1)+(z_1-1)^3+2\left(1-{z_1\over 3}\right)^3\right)\\
&\hspace{-0.2cm}>&\hspace{-0.2cm}{1\over 6}\left({5\over 9}-{4\over 9}\sqrt{1\over {10}}\right). \hspace{7.7cm}(4.59)
\end{eqnarray*}

When $z_2\le y_1-2\le -0.5$, let $P_9$ be the polytope defined in Subcase 3.8 of Lemma 4.4 and define
$$P_{13}=P_1\cap \left\{ (x,y,z):\ 1\le x\le 2,\ y_1-1\le y\le 1, \ z_2+1\le z\le z_1-1\right\}.$$
It can be verified that ${\rm int}(P_8)$, ${\rm int}(P_9)$ and ${\rm int}(P_{13})$ are pairwise disjoint, and
$$P_i\subseteq D'_1(C, X'),\qquad i=8,\ 9,\ 13.$$
Then, by determining the vertices of $P_8$ and $P_{13}$ and routine computations one can deduce
\begin{eqnarray*}
{\rm vol}(P_8)&\hspace{-0.2cm}\ge &\hspace{-0.2cm} {\rm vol}(P_8)\bigg|_{y_2={1\over 2}}={1\over 6}\left({1\over 2}\left(z_2+{3\over 2}\right)^3-2(z_2+1)^3+{1\over 2}\left(z_2+{1\over 2}\right)^3\right)\\
&\hspace{-0.2cm}=&\hspace{-0.2cm} {1\over 6}\left({3\over 4}(z_2+1)-(z_2+1)^3\right),\\
{\rm vol}(P_9)&\hspace{-0.2cm}\ge &\hspace{-0.2cm}{1\over 6}\left(1-(2-y_1)^3-(z_2+1)^3-{1\over 4}(y_1-z_2-2)^3-{3\over 2}(1-y_1-z_2)^2(y_1-1)\right)\\
&\hspace{-0.2cm}\ge &\hspace{-0.2cm}{1\over 6}\left(1-(2-y_1)^3-(y_1-1)^3-{3\over 2}(2-y_1)^2(y_1-1)\right),\\
{\rm vol}(P_{13})&\hspace{-0.2cm}\ge &\hspace{-0.2cm}\left({{4z_2+5}\over 8}+{{4z_2+y_1+3}\over 4}\left({3\over 2}-y_1\right)\right) (z_1-z_2-2)\ge {{4z_2+5}\over 8}(z_1-z_2-2),
\end{eqnarray*}
where $P_{13}$ can be illustrated by a figure similar to Figure 8, and therefore
\begin{eqnarray*}
\hspace{1.5cm}{\rm vol}(D'_1(C, X'))&\hspace{-0.2cm}\ge &\hspace{-0.2cm}{\rm vol}(T_1\cup T_2\cup T_3\cup T_4 \cup P_8\cup P_9\cup P_{13})\\
&\hspace{-0.2cm}\ge &\hspace{-0.2cm}{1\over 6}\left(1+2(x_1-1)^3+(z_1-1)^3-(2-y_1)^3-{3\over 2}(2-y_1)^2(y_1-1)\right.\\
&\hspace{-0.2cm}&\hspace{-0.2cm}\left.+{3\over 4}(z_2+1)-(z_2-1)^3 +{3\over 4}(4z_2+5)(z_1-z_2-2)\right).\\
&\hspace{-0.2cm}\ge &\hspace{-0.2cm}{1\over 6}\left(1+2(x_1-1)^3+(z_1-1)^3-(2-y_1)^3-{3\over 2}(2-y_1)^2(y_1-1)\right.\\
&\hspace{-0.2cm}&\hspace{-0.2cm}\left.+{3\over 4}(y_1-1)-(y_1-1)^3 +{3\over 4}(z_1-y_1)\right).\\
&\hspace{-0.2cm}\ge &\hspace{-0.2cm} {1\over 6}\left( 1+2(x_1-1)^3 +{3\over 4}(z_1-1)+(z_1-1)^3-(y_1-1)^3\right.\\
&\hspace{-0.2cm}&\hspace{-0.2cm}\left. -(2-y_1)^3 -{3\over 2}(2-y_1)^2(y_1-1)\right)\\
&\hspace{-0.2cm}= &\hspace{-0.2cm} {1\over 6}\left( 2(x_1-1)^3 +{3\over 4}(z_1-1)+(z_1-1)^3+{3\over 2}(y_1-1)-{3\over 2}(y_1-1)^3\right)\\
&\hspace{-0.2cm}> &\hspace{-0.2cm}{1\over 6}\left({5\over 9}-{4\over 9}\sqrt{1\over {10}}\right). \hspace{8.5cm}(4.60)
\end{eqnarray*}

As a conclusion of (4.58), (4.59) and (4.60), in this subcase we have
$${\rm vol}(D'_1(C, X'))> {1\over 6} \left({5\over 9}-{4\over 9}\sqrt{1\over {10}}\right).$$

\bigskip
As a conclusion of Subcases 2.1--2.10, when ${\rm int}(F_0)\cap J\not= \emptyset$ we get
$${\rm vol}(D'_1(C,X'))> {1\over 6} \left({5\over 9}-{4\over 9}\sqrt{1\over {10}}\right).$$
Thus, together with Case 1, Lemma 4.5 is proved. \hfill{$\Box $}

\bigskip
Clearly, Theorem 4.1 follows from Lemmas 4.1-4.5.

\bigskip
As a consequence of Lemma 2.1, Lemma 3.2, Theorem 4.1 and (1.2) we get the following upper bounds for $\delta^t(C)$
and $\delta^t(T)$.

\medskip
\noindent
{\bf Theorem 4.2.}
$$\delta^t(C)\le {{90\sqrt{10}}\over {95\sqrt{10}-4}}\approx 0.9601527\cdots\qquad and\qquad
\delta^t(T)\le {{36\sqrt{10}}\over {95\sqrt{10}-4}}\approx 0.3840610\cdots .$$

\vspace{0.5cm}\noindent
{\bf Acknowledgements.} This work is supported by 973 Programs 2013CB834201 and 2011CB302401, the National Natural Science Foundation of China (No. 11071003), and the Chang Jiang Scholars Program of China. For helpful comments, the author is grateful to Prof. P.M. Gruber and Prof. J.C. Lagarias.

\bibliographystyle{amsplain}

\end{document}